\documentclass[11pt, reqno]{amsart}
\usepackage{amsmath,amsthm,amssymb,graphicx}

\theoremstyle{plain}
\newtheorem{thm}{Theorem}[section]
\newtheorem{prop}[thm]{Proposition}
\newtheorem{lem}[thm]{Lemma}
\newtheorem{cor}[thm]{Corollary}

\theoremstyle{definition}
\newtheorem{rem}[thm]{Remark}
\newtheorem{defn}[thm]{Definition}
\newtheorem{eg}[thm]{Example}
\numberwithin{equation}{section}

\newcommand{\bthm}{\begin{thm}}
\newcommand{\ethm}{\end{thm}}
\newcommand{\bprop}{\begin{prop}}
\newcommand{\eprop}{\end{prop}}
\newcommand{\bcor}{\begin{cor}}

\newcommand{\ecor}{\end{cor}}
\newcommand{\blem}{\begin{lem}}
\newcommand{\elem}{\end{lem}}
\newcommand{\bproof}{\begin{proof}}
\newcommand{\eproof}{\end{proof}}
\newcommand{\bca}{\begin{cases}}
\newcommand{\eca}{\end{cases}}
\newcommand{\brem}{\begin{rem}}
\newcommand{\erem}{\end{rem}}
\newcommand{\bpm}{\begin{pmatrix}}
\newcommand{\epm}{\end{pmatrix}}
\newcommand{\bdefn}{\begin{defn}}
\newcommand{\edefn}{\end{defn}}
\newcommand{\bsub}{\begin{subtitle}}
\newcommand{\esub}{\end{subtitle}}
\newcommand{\ben}{\begin{enumerate}}
\newcommand{\een}{\end{enumerate}}
\newcommand{\beg}{\begin{eg}}
\newcommand{\eeg}{\end{eg}}

\def\ss{\smallskip}
\def\ms{\medskip}
\def\bs{\bigskip}
\def\ni{\noindent}

\def\p{\partial}

\def\Im{{\rm Im\/}}

\def\Res{{\rm Res\/}}

\def\Gr{{\rm Gr\/}}
\def\I{{\rm I\/}}
\def\rk{{\rm rk\/}}
\def\ti{\tilde}

\def \a {\alpha}
\def \b {\beta}

\def \e {\epsilon}

\def \l {\lambda}

\def \n {\,\vert\,}
\def \N {\,\Vert\,}
\def \o {\theta}

\def\W{\Omega}

\def\R{\mathbb{R} }
\def\C{\mathbb{C}}

\def\N{\mathbb{N}}

\def\ce{{\mathcal {E}}}

\def\cl{{\mathcal{L}}}
\def\cm{{\mathcal {M}}}

\def\co{{\mathcal {O}}}

\def\cs{{\mathcal {S}}}

\def\tr{{\rm tr}}

\begin{document}


\title[BT and Ward Solitons]
{B\"{a}cklund transformations, Ward solitons,\\
and unitons}
\author{Bo Dai$^*$}\thanks{$^*$Research supported in part by the
AMS Fan Fund\/}
\address{Partner Group of MPI at AMSS\\ Institute of Mathematics\\
Chinese Academy of Sciences\\ Beijing 100080\\ P.~R.~China}
\email{daibo@mail.amss.ac.cn}
\author{ Chuu-Lian Terng$^\dag$}\thanks{$^\dag$Research supported
in  part by NSF Grant DMS-0306446}
\address{
Northeastern University, Boston, MA 02115 and University of
California at Irvine, Irvine, CA 92697} \email{terng@neu.edu}

\ms \hskip 3in 

\begin{abstract}
The Ward equation, also called the modified $2+1$ chiral model, is
obtained by a dimension reduction and a gauge fixing from the
self-dual Yang-Mills field equation on $\R^{2,2}$.  It has a Lax
pair and is an integrable system.  Ward constructed solitons whose
extended solutions have distinct simple poles.  He also used a
limiting method to construct $2$-solitons whose extended solutions
have a double pole. Ioannidou and Zakrzewski, and Anand
constructed more soliton solutions whose extended solutions have a
double or triple pole.  Some of the main results of this paper
are: (i) We construct algebraic B\"acklund transformations (BTs)
that generate new solutions of the Ward equation from a given one
by an algebraic method. (ii) We use an order $k$ limiting method
and algebraic BTs to construct explicit Ward solitons, whose
extended solutions have arbitrary poles and multiplicities. (iii)
We prove that our construction gives all solitons of the Ward
equation explicitly and the entries of Ward solitons must be
rational functions in $x, y$ and $t$.  (iv) Since stationary Ward
solitons are unitons, our method also gives an explicit
construction of all $k$-unitons from finitely many rational
maps from $\C$ to $\C^n$.
\end{abstract}

\maketitle

\bs \tableofcontents

\section{Introduction}

The $2+1$  {\it chiral model \/} is the Euler-Lagrange equation of
the functional
$$\ce(J)=\int_{\R^3} ||J^{-1}J_x||^2 + ||J^{-1}J_y||^2 - ||J^{-1}J_t||^2
\ dx\,dy\,dt,$$ where $||\xi||^2=-\tr(\xi^2)$, $x,y,t$ are the
standard space-time variables, and $J$ is a map from the Lorentz
space $\R^{2,1}$ to the Lie group $SU(n)$.  In other words, $J$ is
a solution of
\begin{equation}\label{aa}
(J^{-1}J_t)_t- (J^{-1}J_x)_x - (J^{-1}J_y)_y =0.
\end{equation}

The {\it Ward equation\/} (or the {\it modified $2+1$ chiral
model\/}) is the following equation for $J:\R^{2,1}\to SU(n)$:
\begin{equation}\label{ward}
(J^{-1}J_t)_t -(J^{-1}J_x)_x -(J^{-1}J_y)_y -[J^{-1}J_t,
J^{-1}J_y]=0.
\end{equation}
This equation is obtained from a dimension reduction and a gauge
fixing of the self-dual Yang-Mills equation on $\R^{2,2}$ (cf.
\cite{W1988}).  We call a solution of the Ward equation a {\it
Ward map\/}.

A Ward map that is independent of $t$ is a harmonic map from
$\R^2$ to $SU(n)$. If the harmonic map has finite energy, then it
extends to a harmonic map from $S^2$ to $SU(n)$. Such harmonic
maps were called {\it unitons\/}, and were studied by Uhlenbeck in
\cite{U1989}, Wood in \cite{Woo1989}, Burstall-Guest in
\cite{BurGue97} and others.

The Ward equation has a {\it Lax pair\/}, i.e., it can be written
as the compatibility condition for a system of linear equations
involving a spectral parameter $\lambda\in\C$. We explain this
next.  Let
\begin{equation} \label{ai}
u=\frac12 (t+y), \quad v=\frac12 (t-y).
\end{equation}
Given smooth maps $A, B:\R^{2,1}\to su(n)$, consider the following
linear system for $\psi:\R^{2,1}\times \C\to GL(n,\C)$:
\begin{equation} \label{eq:auxiliary}
\left\{
\begin{array}{l}
(\lambda\partial_x -\partial_u)\psi =A\psi, \\
(\lambda\partial_v -\partial_x)\psi =B\psi.
\end{array} \right.
\end{equation}
  System
\eqref{eq:auxiliary} is overdetermined. Its compatibility
condition is
$$\left[\l \p_x -\p_u-A, \ \l\p_v -\p_x -B\right]=0.$$
Equate the coefficient of $\l^j$ in the above equation to get
\begin{equation} \left\{
\begin{array}{l}
B_x =A_v, \\
A_x -B_u -[A,B] =0.
\end{array} \right. \label{eq:integrability}
\end{equation}
Suppose $\psi:\R^{2,1}\times \W\to GL(n,\C)$ is a smooth solution
of \eqref{eq:auxiliary} and satisfies the  {\it $U(n)$-reality
condition\/} in $\l$:
\begin{equation}
\psi (x,u,v,\bar{\lambda})^* \psi (x,u,v,\lambda)={\rm I\/},
\label{eq:reality}
\end{equation}
(i.e., $\psi^*=\bar{\psi}^T$), where $\W$ is an open subset of $0$
in $\C$. Let
$$J(x,u,v)=\psi
(x,u,v,0)^{-1}.$$ Then $$A=J^{-1} J_u, \quad B=J^{-1} J_x. $$ Thus
the compatibility condition (\ref{eq:integrability}) implies that
\begin{equation} \partial_v (J^{-1} J_u) =\partial_x (J^{-1} J_x).
\label{eq:xuv} \end{equation} Change back to the standard
variables $(x,y,t)$ to see that $J$ is a solution to the Ward
equation.

A solution $\psi$ to the linear system (\ref{eq:auxiliary}) that
satisfies the $U(n)$-reality condition is called an {\em extended
solution of the Ward equation} or {\it extended Ward map\/}, and
$J=\psi(\cdots, 0)^{-1}$ is the corresponding Ward map.  The
reality condition for $\psi$ implies that $J$ is unitary.  In
other words,  if we find a $\psi(x,y,t,\l)$ so that $\psi$
satisfies the $U(n)$-reality condition and
$(\l\psi_x-\psi_u)\psi^{-1}$ and $(\l\psi_v-\psi_x)\psi^{-1}$ are
independent of $\l$, then $J(x,y,t)=\psi(x,y,t, 0)^{-1}$ is a Ward
map.

The Ward equation has an infinite number of conservation laws
\cite{IW1995}. In particular, the energy functional
$$E(J)= \frac12 \iint_{\{ t={\rm const}\} } \Vert J^{-1} J_t \Vert^2
+\Vert J^{-1} J_x \Vert^2 +\Vert J^{-1} J_y \Vert^2 \ dxdy,
$$
is a conserved quantity. To ensure finite energy, Ward imposed the
following boundary condition
\begin{equation} \label{bc1}
J=J_0 +J_1(\theta) r^{-1} +O(r^{-2}) \quad {\rm as\/}\ r \to
\infty,
\end{equation} where $x+iy =r
e^{i\theta}$, $J_0$ is a constant matrix, and $J_1$ is independent
of $t$.

A Ward map $J$ is called a {\it Ward soliton\/} if  $J$ \ben \item
has finite energy on $\R^2$, or equivalently satisfies the
boundary condition \eqref{bc1}, \item has an extended solution
$\psi$ such that $\psi(x,y,t,\l)$ is rational in $\l$ and
$$\lim_{|\l|\to \infty}\psi(x,y,t,\l)= {\rm I\/}$$
for all $(x,y,t)$. \een

If $\psi$ is an extended solution of the Ward equation with poles
at $\l= z_1, \ldots, z_r$ of multiplicities $n_1, \ldots, n_r$
respectively, then $(z_1, \ldots, z_r, n_1, \ldots, n_r)$ is
called the {\it pole data\/} of $\psi$ and $\sum_{j=1}^r n_j$ is
called the {\it degree\/} of $\psi$. A Ward soliton $J$ is called
a $k$-soliton if $k$ is the minimum of
$$\{\deg(\psi)\n \psi \ {\rm is\ an\  extended\ solution\  of \ }J\}.$$

Let $z\in \C\setminus \R$, $\pi$ a Hermitian projection of $\C^n$,
and
$$g_{z,\pi}(\l) =\pi+ \frac{\l-\bar z}{\l-z}\pi^\perp= \I +
\frac{z-\bar z}{\l-z} \ \pi^\perp.$$ A direct computation shows
that $g_{z,\pi}$ satisfies the $U(n)$-reality condition
\eqref{eq:reality}.  Such $g_{z,\pi}$ is called a {\it simple
element\/}.

Let $z\in \C\setminus \R$ be a fixed constant, $\cm_{n\times k}^0$
the space of rank $k$ complex $n\times k$ matrices, and
$V=(v_{ij}):\C\to \cm^0_{n\times k}$ a meromorphic map. Let $\pi$
denote the map from $\R^{2,1}$ to the space of rank $k$ Hermitian
projections of $\C^n$ such that $\Im(\pi(x,y,t))$ is the complex
linear subspace of $\C^n$ spanned by columns of
$$V(x+zu+z^{-1}v),$$
and $\pi^\perp= \I-\pi$. Ward (cf. \cite{W1988}) noted that
\begin{equation} \label{bd}
g_{z,\pi(x,y,t)}=\I +\frac{ z-\bar z}{\l-z}\pi^\perp(x,y,t),
\end{equation}
is an extended solution, i.e., a solution of \eqref{eq:auxiliary}
with
$$A= (\bar z-z)\pi_x, \quad
B=(\bar z-z)\pi_v.$$ The associated Ward map is
$$J_{z,V}(x,y,t) = g_{z, V}(x,y,t,0)^{-1}=\pi(x,y,t)+ \frac{z}{\bar
z} \pi^\perp(x,y,t).$$ Ward proved that $J_{z,V}$ satisfies the
boundary condition (\ref{bc1}) if and only if each $v_{ij}$ is a
rational function (cf. \cite{W1988}). Hence $J_{z,V}$ is a Ward
$1$-soliton if each entry of $V$ is a rational function. Note
$J_{i,V}$ is a stationary Ward map, i.e.,
 a harmonic map from $\C$ to $U(n)$.

There are several methods for constructing exact Ward
multi-solitons: Ward used the method of Riemann-Hilbert problem
with zeros in \cite{W1988} to construct $k$-soliton solutions
whose extended solutions have $k$ simple poles.  Such solutions
have trivial scattering in the sense that the $k$ one-solitons
preserve their travelling directions and shapes after the
interaction. Taking the limit of an extended $2$-soliton with
poles at $i+\e$ and $i-\e$ as $\e\to 0$,  Ward and Ioannidou found
extended $2$-solitons with a double pole at $\l=i$ (cf.
\cite{W1995,I1996,IZ1998}). Ioannidou also constructed some
extended $3$-solitons with a triple pole at $\l=i$. These limiting
solutions have non-trivial scattering,  i.e., the travelling
directions of interacting localized lumps change after the
interaction. For example, Ioannidou give examples of extended
$2$-solitons with a double pole at $\l=i$  and with scattering
angle $\pi/k$.  Anand constructed more solitons with non-trivial
scattering in \cite{ A1997, Ana98}. Ioannidou and Zakrzewski
generalized Uhlenbeck's method of adding unitons for harmonic map
equation to Ward equation in \cite{IZ1998} by writing down an
analytic B\"acklund transformation. Vilarroel, Fokas and Ioannidou
studied the inverse scattering of the Ward equation in
\cite{V1990,FokIoa98}. Zhou gave Darboux transformations in
\cite{Zh93}.

The standard analytic B\"acklund transformations (BT)  goes as
follows:  Given an extended solution $\psi$ of
\eqref{eq:auxiliary}, if we want to find a projection map $\ti
\pi$ so that $\psi_1= g_{z, \ti \pi}\psi$ is again an extended
solution, then the condition that $\psi_1$ satisfies
\eqref{eq:auxiliary} for some $\ti A(x,y,t)$ and $\ti B(x,y,t)$ is
equivalent to the condition that $\ti\pi$ is a solution of the
following  system of first order partial differential equations:
\begin{equation*}
\bca
\ti\pi^\perp(z\ti\pi_x-\ti\pi_u-A\ti\pi)=0,&\\
\ti\pi^\perp(z\ti \pi_v-\ti\pi_x-B\ti\pi)=0, \eca \hskip 1in ({\rm
BT\/}_{z,\psi})
\end{equation*}
where $A=(\l\psi_x-\psi_u)\psi^{-1}$ and
$B=(\l\psi_v-\psi_x)\psi^{-1}$.
 A solution of ${\rm BT\/}_{z,\psi}$ gives rise to an explicit
extended Ward map with one extra pole at $\l= z$. Although this
first order PDE is solvable,  general solutions have not been
fully understood.  One result of this paper is an explicit
construction of all solutions of  BT$_{z,\psi}$ when  $\psi$ is an
extended Ward soliton.

Another result of this paper is to construct  an algebraic BT for
the Ward equation. This is a transformation that generates a new
extended solution $\psi_1$ by an algebraic formula in terms of a
given extended solution $\psi$ and an extended $1$-soliton $g_{z,
\pi}$.  In fact, if $\psi$ is holomorphic and non-degenerate at
$\l= z$, then
$$\psi_1(x,y,t,\l)= g_{z,\ti\pi(x,y,t)}\psi(x,y,t,\l)= \left({\rm
I\/} + \frac{z-\bar z}{\l-z}\ \ti\pi^\perp (x,y,t)
\right)\psi(x,y,t,\l)$$ is also an extended solution of the Ward
equation, where $\ti\pi(x,y,t)$ is the Hermitian projection onto
$\psi(x,y,t,z){\rm Im\/}(\pi(x,y,t))$. In other words, $\ti \pi$
is a solution of ${\rm BT\/}_{z,\psi}$.  Note that the algebraic
BT only works if the given extended solution $\psi$ is holomorphic
and non-degenerate at $\l=z$.  In this case, the new extended
solution $\psi_1$ has one more pole at $\l=z$ than $\psi$. We
apply algebraic BTs repeatedly to an extended $1$-soliton to get
Ward's multi-solitons, whose extended solutions have distinct
poles. We use algebraic BTs $k$ times and a delicate limiting
method to construct multi-solitons, whose extended solutions have
general pole data $(z_1, \ldots, z_r, n_1, \ldots, n_r)$.

There are also analytic and algebraic BTs for harmonic maps from
$\R^2$ to $U(n)$ (\cite{U1989, BerGue91}).  But the algebraic BT
of a finite energy harmonic map has infinite energy.  Hence we
cannot produce new harmonic maps on $S^2$ using algebraic BTs.
Although Uhlenbeck's adding uniton method  can be viewed as the
limiting case of algebraic BTs as the pole goes to $i$, the limit
of these BTs of a harmonic map $s$ gives the same $s$ (for more
detail, cf. \cite{U1989, BerGue91}). However, if we apply
algebraic BTs of the Ward equation with pole at $i+\e$ to a
$1$-uniton and choose the projection $\pi_\e$ of the Ward
$1$-soliton $g_{i+\e,\pi_\e}$ carefully, then  as $\e\to 0$ the
limiting solution can be a $2$-uniton.  We show in this paper that
this limiting method for the Ward equation can produce all unitons
into $U(n)$. In fact, we give an explicit construction of
$k$-unitons from $k$ rational maps from $\C$ to $\C^n$. Our
construction of unitons is different from the ones given by Wood
in \cite{Woo1989} and by Burstall-Guest in \cite{BurGue97}.

This paper is organized as follows:  We give a quick review of
unitons and Ward $1$-solitons in section 2, give algebraic
B\"acklund transformations for the Ward equation in section 3.
Uhlenbeck proved that a rational map $f:S^2\to GL(n,\C)$
satisfying the $U(n)$-reality condition \eqref{eq:reality} and
$f(\infty)=\I$ can be factored as a product of simple elements.
But such factorization in general is not unique.   We give a
refinement of this factorization so that it is unique in section
4. We apply B\"acklund transformations and a careful limiting
method to construct Ward solitons that satisfy the boundary
condition \eqref{bc1} and their extended solutions have pole data
$(z,k)$ in section 5.    We construct multi-solitons whose
extended solutions have pole data $(z_1, \ldots, z_r, n_1, \ldots,
n_r)$ in section 6.  We show in section 7 that the first equation
of ${\rm BT\/}_{z,\psi}$ defines a natural complex structure on
the trivial bundle $S^2 \times \C^n$ over $S^2$, and  a solution
of ${\rm BT\/}_{z,\psi}$ corresponds to a holomorphic subbundle of
the trivial bundle  that  satisfies certain first order PDE
constraint. In section 8, we  use the holomorphic vector bundle
formulation of section 7 to  prove that algebraic BTs and the
limiting method of section 5 produce all solutions of ${\rm
BT}_{z,\psi}$ for any extended Ward-soliton $\psi$, hence we can
construct all Ward solitons explicitly.   In section 9, we give an
explicit construction of all unitons using the limiting method of
section 5.   

The graphics of Ward solitons indicate that a
Ward soliton with polo data $(z_1, \ldots, z_k, n_1, \ldots, n_k)$
is the interaction of $k$ Ward solitons with pole data $(z_1,
n_1), \ldots$, $(z_k,n_k)$ respectively and these $k$ solitons
keep their shapes after interaction.  But  the dynamics of
solitons with pole data $(z,k)$ are intriguing,  quite
complicated, and deserve further investigation.
The reader can  play the Quick Time movies for several examples of Ward solitons 
by going to the following website:

\centerline{http://www.math.neu.edu/$\sim$terng/WardSolitonMovies.html.}

The first author would like to thank
the AMS Fan Fund and Northeastern University for sponsoring his
visit to Northeastern University, where the cooperation started.
The second author also thanks Karen Uhlenbeck for many useful
discussions, and thanks MSRI for supporting her visit during the
winter quarter of 2004, where she worked on this paper.

\bs
\section{1-unitons and 1-soliton Ward maps}

A stationary solution of the Ward equation is a harmonic map from
$\R^2$ to $U(n)$.  If in addition it has finite energy then it is
a harmonic map from $S^2$. All such harmonic maps are called
unitons, which are studied by Uhlenbeck \cite{U1989}, Wood
\cite{Woo1989}, Burstall-Guest \cite{BurGue97} and others.

The harmonic map equation is integrable in the sense that there is
an associated linear system with a complex parameter $\xi\in
\C\setminus \{0\}$.  Namely, if $s$ is a harmonic map from $\C$ to
$U(n)$, then the following linear system is compatible:
\begin{equation}\label{ae}
\bca E_z= (1-\xi^{-1}) E  P, &\cr E_{\bar z}=-(1-\xi) E P^*, \eca
\end{equation}
where $P=\frac{1}{2} s^{-1}s_z$.  Note that the compatibility
condition of system \eqref{ae} is the harmonic map equation,
$$P_{\bar z}=- [P, P^*].$$
Conversely, if $E(x, y, \xi)$ is a solution of \eqref{ae} and
satisfies the $U(n)$-reality condition \eqref{eq:reality}, then
$s(x,y)=E(x,y,-1)$ is a harmonic map.  Such $E$ is called an {\it
extended solution\/} of the harmonic map equation. A direct
computation implies that if $E(x,y,\xi)$ is an extended solution
of the harmonic map equation, then
$$\psi(x,y,t,\l)= E\left(x,y, \frac{\l-i}{\l+i}\right)^{-1}$$ is
an extended solution of the Ward equation, i.e., $\psi$ is a
solution of \eqref{eq:auxiliary} and $\psi(x, y, t,
0)^{-1}=s(x,y)$ is a stationary Ward map.

Let  $V=(v_{ij}):\C\to \cm_{n\times k}^0(\C)$ be a rational map,
$\pi$ the projection of $\C^n$ onto the subspace spanned by the
$k$ columns of $V$, and $\pi^\perp=\I-\pi$.  Then
$s=\pi-\pi^\perp$ is a $1$-uniton.  Moreover, all $1$-unitons are
of this form.  The $1$-uniton $\pi-\pi^\perp$ has an extended
solution:
\begin{equation}\label{bb}
E(x,y,\xi)= \pi(x,y) + \xi \pi(x,y)^\perp.
\end{equation}
Uhlenbeck proved in \cite{U1989} that given a harmonic map
$s:S^2\to U(n)$, there exists an extended solution $E(x,y,\xi)$ of
the form
\begin{equation}\label{bs}
E(x,y,\xi)=(\pi_1+\xi\pi_1^\perp)\cdots (\pi_k+\xi\pi_k^\perp),
\end{equation}
where each $\pi_i(x,y)$ is a projection onto some
$k_i$-dimensional linear subspace $V_i(x,y)$ of $\C^n$ and $k\leq
(n-1)$.  Such solutions are called {\it k-unitons\/}.

Substitute $\xi= \frac{\l-i}{\l+i}$ into \eqref{bs} to get an
extended Ward $k$-soliton with pole data $(i,k)$.  In particular,
\begin{equation}\label{ba}
\psi(x,y,t,\l) =E(x,y,\frac{\l+i}{\l-i})= \pi(x,y)
+\frac{\l+i}{\l-i}\ \pi^\perp (x,y) = g_{i, \pi(x,y)}(\l),
\end{equation}
is an extended Ward $1$-soliton with a simple pole at $\l=i$.

    A general extended Ward $1$-soliton \eqref{bd} is obtained by
replacing $i$ by a non-real complex constant $z$ and $x+iy$ by
$w=x+zu+z^{-1}v$.  The associated Ward map is
$$\hat J_{z,V}(x,y,t)=\psi(x,y,t, 0)^{-1} =\pi +\frac{z}{\bar z}
\pi^\perp .$$ Note that $\hat J_{z,V}$ has constant determinant
$(z/ \bar{z})^{n-k}$. So we can normalize it to get a Ward map
into $SU(n)$:
$$J_{z,V}(x,y,t)= \left( \frac{\bar z}{z} \right)^{\frac{n-k}{n}}
\left(\pi +\frac{z}{\bar z} \pi^\perp \right)= \left(
\frac{z}{\bar z} \right)^{\frac{k}{n}} \left(\frac{\bar z}{ z} \pi
+ \pi^\perp \right).$$

    The $1$-soliton  $J_{z,V}$ is a travelling wave. To see this,
write $z=r e^{i\o}$ and compute directly to get
$$w=x+zu+z^{-1}v =(x-v_1 t) +k_1(y-v_2 t) +ik_2(y-v_2 t),$$
where $v_1=-\frac{2r \cos\theta}{1+r^2}$,
$v_2=\frac{1-r^2}{1+r^2}$, and $k_1 +ik_2=(z -z^{-1})/2$. Thus
$J_{z,V}$ is a travelling wave with constant velocity
\begin{equation} \label{ek}
\vec{v} =(-\frac{2r \cos\theta}{1+r^2}, \frac{1-r^2}{1+r^2})
\end{equation}
 on the $xy$-plane.

\beg Ward $1$-solitons
\par

Let $z\in\C \setminus \R$, $f:\C\to \C$ a rational function,
$w=x+zu + z^{-1}v$, $V(w)= \bpm 1\cr f(w)\epm$, and $\pi(x,y,t)$
the projection onto $\C V(w)$.  A direct computation gives
$$J_{z,V}= \frac{1}{|z|(1+|f(w)|^2)} \bpm \bar z+ z
|f(w)|^2 & (\bar z-z) \overline{f(w)}\cr (\bar z-z) f(w) & \bar z
|f(w)|^2 + z\epm.$$

\eeg

\bs
\section{Algebraic B\"{a}cklund transformations (BT)}

 In this section, we give an algebraic BT to construct a
family of explicit solutions from a given  extended solution
$\psi(\l)(x,y,t)=\psi(x,y,t,\l)$ of the Ward equation.

\begin{thm}[Algebraic B\"{a}cklund transformation]
\label{thm:BT} Let $\psi(x,y,t,\l)$ be an extended solution of the
Ward equation, and $J= \psi(\cdots, 0)^{-1}$ the associated Ward
map. Choose $z\in {\mathbb C}\setminus {\mathbb R}$ such that
$\psi$ is holomorphic and non-degenerate at $\lambda =z$. Let
$g_{z,\pi(x,y,t)}(\lambda)$ be an extended $1$-soliton, and
$\ti\pi(x,y,t)$ the Hermitian projection of $\C^n$ onto
$$\psi(x,y,t,z)(\Im(\pi(x,y,t)).$$ Then
\ben \item $\ti\psi (x,y,t,\lambda) =
g_{z,\tilde{\pi}(x,y,t)}(\lambda)
\psi(x,y,t,\lambda)g_{z,\pi(x,y,t)}(\l)^{-1}$ is holomorphic and
non-degenerate at $\l= z, \bar z$, \item $\psi_1= g_{z,\ti\pi}
\psi = \ti\psi g_{z,\pi}$ is a new extended solution to the linear
system (\ref{eq:auxiliary}) with $$(A,B)\to (A +(\bar{z}-z)
\tilde{\pi}_x, B + (\bar{z}-z) \tilde{\pi}_v),$$ and the new Ward
map is
$$J_1(x,y,t) =\left( \frac{z}{\bar{z}} \right)^{k/n}
J(x,y,t)\, \left( \frac{\bar{z}}{z} \tilde{\pi} (x,y,t)
+{\tilde{\pi}}^\perp (x,y,t) \right).$$ \een
\end{thm}

\begin{proof}
 (1) Let $\tilde{\psi}(\lambda) =g_{z,\tilde{\pi}} (\lambda)
\psi(\lambda) g_{z,\pi}(\lambda)^{-1}$. Then residue calculus
implies that $\tilde{\psi}(\l)$ is holomorphic at $\lambda
=z,\bar{z}$. Thus we have two factorizations of $\psi_1
=g_{z,\tilde{\pi}} \psi =\tilde{\psi} g_{z,\pi}$.

(2) It suffices to show that
\begin{equation}
A_1:= (\lambda\partial_x \psi_1 -\partial_u \psi_1) \psi_1^{-1}
\label{expr1} \end{equation} is independent of $\lambda$. Using
$\psi_1 =g_{z,\tilde{\pi}} \psi$, we have
\begin{equation} A_1= (\lambda \partial_x  g_{z,\tilde{\pi}}
-\partial_u  g_{z,\tilde{\pi}}) g_{z,\tilde{\pi}}^{-1} +
g_{z,\tilde{\pi}}(\lambda \partial_x \psi -\partial_u
\psi)\psi^{-1} g_{z,\tilde{\pi}}^{-1}. \label{expr2}
\end{equation} Since $(\lambda \partial_x \psi -\partial_u
\psi)\psi^{-1}$ is constant in $\lambda$ by assumption,
(\ref{expr2}) is holomorphic at $\l\in {\mathbb C} \setminus \{
z,\bar{z}\}$, and has at most a simple pole at $\lambda =\infty$.
But $$\Res_{\lambda=\infty} A_1 =\left. (\partial_x
g_{z,\tilde{\pi}}) g_{z,\tilde{\pi}}^{-1} \right\vert
_{\lambda=\infty} =0$$ as $g_{z,\tilde{\pi}} (\infty) \equiv {\rm
I\/}$. So $A_1$ is holomorphic at $\l\in {\mathbb C}\cup
\{\infty\} \setminus \{ z,\bar{z}\}$. On the other hand, using
$\psi_1 =\tilde{\psi} g_{z,\pi}$, \begin{equation} A_1= (\lambda
\partial_x \tilde{\psi} -\partial_u  \tilde{\psi})
\tilde{\psi}^{-1} + \tilde{\psi} (\lambda \partial_x g_{z,\pi}
-\partial_u g_{z,\pi}) g_{z,\pi}^{-1} \tilde{\psi}^{-1}.
\label{expr3} \end{equation} Since $(\lambda \partial_x g_{z,\pi}
-\partial_u g_{z,\pi}) g_{z,\pi}^{-1}$ is independent of
$\lambda$, $A_1$ is holomorphic at $\lambda =z,\bar{z}$. Thus we
see that $A_1$ is holomorphic on ${\mathbb C}\cup \{\infty\}$,
hence independent of $\lambda$ by Liouville's Theorem. Likewise
$(\lambda\partial_v \psi_1 -\partial_x \psi_1) \psi_1^{-1}$ is
also independent of $\lambda$. The remaining computation is
straightforward.
\end{proof}

 Let
$$\psi_1=g_{z,\pi} \ast \psi, \quad J_1=g_{z,\pi} \ast J$$
denote the algebraic B\"{a}cklund transformation generated by
$g_{z,\pi}$. If we apply BTs repeatedly (with distinct poles) to
an extended $1$-soliton solution, then we obtain Ward
multi-solitons, whose extended solutions have only simple poles.
Such solutions coincide with the ones obtained by Ward
\cite{W1988} using solutions of the Riemann-Hilbert problem.

\beg\label{af}  Ward $2$-solitons with trivial scattering
\par

Let $z_1, z_2$ be two distinct  complex numbers and $z_1\not=
\bar{z}_2$, $f_1, f_2:\C\to \C$ rational functions, and
$\pi_i(x,y,t)$ the projection onto $\C\bpm 1\cr f_i(w_i)\epm$,
where
$$w_i= x+ z_i u + z_i^{-1}v, \quad i=1,2.$$
Then $g_{z_1, \pi_1}$ and $g_{z_2,\pi_2}$ are extended $1$-soliton
solutions of the Ward equation.  Apply B\"{a}cklund transformation
(Theorem \ref{thm:BT}) with $\psi= g_{z_1, \pi_1}$ and $g_{z,\pi}=
g_{z_2, \pi_2}$. Compute directly to see that $\psi(z_2)\left(\bpm
1 \cr f_2(w_2)\epm\right)$ is parallel to
\begin{equation}\label{ah}
\ti v_2= A \bpm 1\\ f_1(w_1)\epm + B\bpm \overline{f_1(w_1)}\\
-1\epm,
\end{equation}
where
$$A=1+\overline{f_1(w_1)} f_2(w_2),
\quad B=\frac{z_2-\bar z_1}{z_2-z_1} (f_1(w_1)-f_2(w_2)). $$ The
new extended solution of the Ward equation is
$$\ti \psi= g_{z_2, \ti \pi_2} g_{z_1, \pi_1},$$
where $\ti \pi_2$ is the projection onto $\C \ti v_2$. The
associated Ward map is
$$J= c(\pi_1+ \frac{z_1}{\bar z_1}\pi_1^\perp) (\ti \pi_2 +
\frac{z_2}{\bar z_2} \ti{\pi}_2^\perp),$$ where $c^2= \frac{\bar
z_1 \bar z_2}{z_1 z_2}$ is a normalizing constant to make
$\det(J)=1$. \eeg

\ss \ni {\bf Remark.}  If $J$ is a Ward map into $SU(n)$, then the
Ward map associated to $\psi_1$ in Theorem \ref{thm:BT} is
$$\hat J_1= J\, (\ti\pi+ \frac{z}{\bar z} \ti\pi^\perp),$$
which is a Ward map into $U(n)$.  But $\det(\hat J_1)= (z/\bar
z)^{n-k}$ is constant.  So $J_1= (\bar z/z)^{\frac{n-k}{k}}\hat
J_1$ is a Ward map into $SU(n)$.  This means that BTs are defined
for both the $SU(n)$ and the $U(n)$ case.

\bs

\section{Minimal factorization\/}

Let $\cs_r(S^2, GL(n))$ denote the group of rational maps
$f:S^2\to GL(n,\C)$ that satisfies the reality condition
$f(\bar\l)^*f(\l)=\I$ and $f(\infty)=\I$. First we recall the
factorization theorem of Uhlenbeck \cite{U1989}

 \bthm \label{bf}  \cite{U1989}
The group  $\cs_r(S^2,GL(n))$ is generated by the set of all
simple elements, i.e., every $f\in \cs_r(S^2,GL(n))$ can be
factored as a product of simple elements,
$$f=g_{z_1, \pi_1}\cdots g_{z_k,\pi_k},$$
for some $z_1,\ldots, z_k$ and Hermitian projections $\pi_1,
\ldots, \pi_k$. \ethm

However, the above factorization is not unique.  For example, if
$\Im \pi_1$ is orthogonal to $\Im\pi_2$ then
$g_{z,\pi_1^\perp}g_{z,\pi_2^\perp}= g_{z,\pi^\perp}$, where $\pi$
is the projection onto $\Im\pi_1\oplus \Im\pi_2$. Moreover, if
$z_1\not=z_2, \bar z_2$, then $g_{z_1,\pi_1}g_{z_2\pi_2}$ can be
written as $g_{z_2,\tau_2}g_{z_1,\tau_1}$ for some projections
$\tau_1, \tau_2$. This is the permutability formula for simple
elements given in  Theorem 6.2 of \cite{TU2000}, which can be
reformulated as follows:

\bthm\label{bg}\cite{TU2000} Suppose $z_1\not= z_2, \bar z_2$, and
$\pi_1, \pi_2$ are Hermitian projections of $\C^n$.  Let
$\ti\pi_1$ be the projection onto $g_{z_2,\pi_2}(z_1)(\Im\pi_1)$,
and $\ti\pi_2$ the projection onto $g_{z_1,\pi_1}(z_2)(\Im\pi_2)$.
Then
$$g_{z_1,\ti\pi_1}g_{z_2,\pi_2}= g_{z_2, \ti\pi_2}g_{z_1, \pi_1}.$$
Conversely, if $\tau_i$ are projections so that
$g_{z_1,\tau_1}g_{z_2,\pi_2}= g_{z_2,\tau_2}g_{z_1,\pi_1}$, then
$\tau_i=\ti\pi_i$ for $i=1, 2$. \ethm

Recall that $g_{z,\pi}\ast \psi$ is the algebraic BT of $\psi$
generated by the $1$-soliton $g_{z,\pi}$. As a consequence of
Theorem \ref{bg} we have

\bcor
 If $g_{z_i, \pi_i}$ are extended Ward $1$-solitons, then
$g_{z_2,\pi_2}\ast g_{z_1, \pi_1}= g_{z_1, \pi_1}\ast g_{z_2,
\pi_2}$.
 \ecor

Note that the proof of Theorem~\ref{thm:BT} (1) gives a more
general permutability formula:

\bprop \label{bm} Suppose $f\in \cs_r(S^2,GL(n))$ is holomorphic
and non-degenerate at $\l=z$ and  $g_{z,\pi}$ is a simple element.
Let $\ti\pi$ be the projection onto $f(z)(\Im\pi)$.  Then \ben
\item $\ti f= g_{z,\ti\pi} f  g_{z,\pi}^{-1}$ is holomorphic at
$\l=z,\bar z$, \item $\ti f g_{z,\pi}= g_{z,\ti\pi} f$. \een
\eprop

It follows from Theorem \ref{bf} and
permutability formula \ref{bg} that we have

\bcor
 Let $f\in \cs_r(S^2,GL(n))$, and $C(f)$ the set of poles of $f$.
Suppose $C(f)\subset \C_+=\{r+is\n s>0\}$, $C_1, C_2$ are proper
disjoint subsets of $C(f)$, and $C(f)= C_1\cup C_2$.  Then there
exist unique $f_1, f_2\in \cs_r(S^2,GL(n))$ so that  $\psi=f_1f_2$
and $C(f_i)=C_i$ for $i=1, 2$. \ecor

\bcor \label{bj} Let $f\in \cs_r(S^2,GL(n))$ with pole data $(z_1,
\ldots, z_k, n_1, \ldots, n_k)$ and $z_i\in \C_+$ for $1\leq i\leq
k$. Then: \ben \item  There exist unique $f_i\in \cs_r(S^2,GL(n))$
with pole data $(z_i, n_i)$ and $F_i\in \cs_r(S^2,GL(n))$ that is
holomorphic and non-degenerate at $z_i, \bar z_i$ such that $f=F_i
f_i$ for each $i$.
 \item There exist $g_j\in
\cs_r(S^2,GL(n))$ with pole data $(z_j,n_j)$ for $1\leq j\leq k$
so that  $f=g_1\cdots g_k$. \een \ecor

Next we give a refinement of the factorization for elements in
$\cs_r(S^2,GL(n))$ whose pole data is $(z,k)$.  First we need a
Lemma.

\begin{lem} \label{ct}
Let $\pi_1$ and $\pi_2$ be two Hermitian projections of $\C^n$
onto $V_1$, $V_2$ respectively. \ben \item
 If $V_1\perp V_2$, then
\begin{equation} \label{ca} g_{z,\pi_2}
g_{z,\pi_1}=\frac{\l-\bar z}{\l-z}\, g_{z,\tau},\end{equation}
where $\tau =\pi_2 +\pi_1$ is the projection onto $V_2\oplus V_1$.
\item Suppose  $V_2^1 := V_2 \cap V_1^\perp \not= 0$.  Let
$\tau_2$ and $\tau_1$ be the  projections onto
$V_2\cap(V_2^1)^\perp$ and $V_1\oplus V_2^1$ respectively.  Then
$\Im\tau_2\cap \Im\tau_1^\perp=0$ and
\begin{equation} \label{cb} g_{z,\pi_2} g_{z,\pi_1}= g_{z,\tau_2}
g_{z,\tau_1},
\end{equation}
\een
\end{lem}

\begin{proof}
A direct computation gives \eqref{ca} and \eqref{cb}. Compute
directly to see
\begin{align*}
(\Im\tau_2)\cap (\Im\tau_1^\perp) &= (V_2\cap (V_2^1)^\perp)
\cap(V_1\oplus V_2^1)^\perp \\ &= (V_2\cap (V_2^1)^\perp) \cap
(V_1^\perp \cap (V_2^1)^\perp)\\ &= (V_2\cap V_1^\perp) \cap
(V_2^1)^\perp =V_2^1 \cap (V_2^1)^\perp =0.
\end{align*}
\end{proof}

\bprop\label{cr} Suppose $\pi_1, \ldots, \pi_k$ are Hermitian
projections of $\C^n$ and $\Im\pi_j\cap \Im\pi_{j-1}^\perp=0$ for
all $2\leq j\leq k$.  Let $n_j=$ the rank of $\pi_j$.  Then \ben
\item $n_1\geq n_2\geq\cdots \geq n_k$, \item ${\rm
Ker\/}(\pi_k^\perp\cdots \pi_1^\perp) = \Im\pi_1$, \item
$\dim(\Im(\pi_k^\perp\cdots \pi_1^\perp)) = n-n_1$. \een \eprop

\begin{proof}
Denote $V_j=\Im \pi_j$, $1\le j\le k$. The kernel of
$\pi_j^\perp:V_{j-1}^\perp\to V_j^\perp$ is $V_j\cap
V_{j-1}^\perp=0$. So $\pi_j^\perp$ is injective on
$V_{j-1}^\perp$.
\end{proof}

\bdefn  Suppose $\phi\in \cs_r(S^2,GL(n))$ has pole data $(z,k)$.
A factorization of $\phi$  is called {\it minimal\/} if
$$\phi =\left( \frac{\l-\bar z}{\l-z} \right)^{k-l}
g_{z,\pi_l}\cdots g_{z,\pi_1}$$ with $\pi_j\not= 0, \I$, and
$\Im\pi_j \cap \Im\pi_{j-1}^\perp = 0$ for $j=2,\cdots l$. \edefn

\bthm \label{cc} If $\phi\in \cs_r(S^2,GL(n))$ has pole data
$(z,k)$, then $\phi$  has a unique minimal factorization. \ethm

\begin{proof} By Uhlenbeck's factorization Theorem \ref{bf},
we can factor
$$\phi_k =g_{z,\pi_k}\cdots g_{z,\pi_1}.$$
We first prove the existence of minimal factorization by induction
on $k$. For $k=1, 2$, the Theorem is true.  Suppose the Theorem is
true for $k-1$. Induction hypothesis implies that
$$g_{z,\pi_{k-1}}\cdots g_{z,\pi_1}=\left(\frac{\l-\bar
z}{\l-z}\right)^m g_{z,\tau_{k-1-m}}\cdots g_{z,\tau_1}$$ so that
the right hand side is a minimal factorization.  If $m\geq 1$,
then by induction hypothesis $g_{z,\pi_k}g_{z,\tau_{k-1-m}}\cdots
g_{z,\tau_1}$ has a minimal factorization.  So does $\phi_k$.  If
$m=0$, then there are two cases: \ben \item  If $\Im\pi_k\cap
\Im\tau_{k-1}^\perp=0$, then $g_{z,\pi_k}g_{z,\tau_{k-1}}\cdots
g_{z,\tau_1}$ is a minimal factorization for $\phi_k$. \item
 If $V:=\Im\pi_k\cap\Im\tau_{k-1}^\perp\not=0$, then by Lemma
\ref{ct} (2) we can write
 $$g_{z,\pi_k}g_{z,\tau_{k-1}}= g_{z,\ti\pi_k} g_{z,\ti\tau_{k-1}}$$
 such that $\Im\ti\pi_k\cap\Im\ti \tau_{k-1}^\perp=0$, where
$\Im\ti\pi_k=\Im\pi_k\cap V^\perp$, and $\Im\ti\tau_{k-1}=
\Im\tau_{k-1}\oplus V$. Since $V\not=0$,
$\rk(\ti\pi_k)<\rk(\pi_k)$. By induction hypothesis,
$g_{z,\ti\tau_{k-1}}g_{z,\tau_{k-2}}\cdots g_{z,\tau_1}$ has a
minimal factorization
$$g_{z, \hat \tau_{k-1}}\cdots g_{z,\hat\tau_1}.$$
If $\Im\ti\pi_k\cap\Im\hat \tau_{k-1}^\perp \not=0$, then we use
Lemma \ref{ct} again to reduce the rank of $\ti\pi_k$.  So after
finitely many times, we can obtain a minimal factorization of
$\phi_k$. \een

Next we use induction on $k$ to prove the uniqueness of minimal
factorization. The case $k=1$ is obvious. Suppose all $\phi$ with
pole data $(z,k)$ and $k<K$ have unique minimal factorizations.
Consider two minimal factorizations
$$\left( \frac{\l-\bar z}{\l-z}
\right)^{K-l}g_{z,\pi_l}\cdots g_{z,\pi_1} =\left( \frac{\l-\bar
z}{\l-z} \right)^{K-m}g_{z,\tau_m}\cdots g_{z,\tau_1}.$$ Compare
the coefficients of $\frac{1}{(\l-z)^K}$ to get
$$(z-\bar z)^{K-l} \pi_l^\perp\cdots \pi_1^\perp= (z-\bar z)^{K-m}
\tau_m^\perp\cdots \tau_1^\perp.$$ By Proposition \ref{cr} (2),
the kernel of the left hand side and the right hand side operators
are $\Im\pi_1$ and $\Im\tau_1$ respectively.  Hence $\pi_1=
\tau_1$. Then induction hypothesis gives the uniqueness.
\end{proof}

\bs

\section{Ward solitons with pole data $(z,k)$}

Ward noted that the limit of extended $2$-soliton solutions with
poles at $i+\e$ and $i-\e$ as $\e\to 0$ gives time dependent Ward
maps, and are $2$-solitons with non-trivial scattering (cf.
\cite{W1995,I1996,IZ1998}).  In this section, we use a systematic
limiting method and algebraic BTs to construct extended Ward
$k$-solitons with pole data $(z,k)$.

First we give the Example of Ward: \beg\label{ag}  Ward
$2$-solitons with non-trivial scattering.
\par
Let $\a\in \C \setminus\R$, and $f, g$ two rational functions on
$\C$. Choose $z_1= \a+ \e$, $z_2=\a-\e$, $f_1= f+\e g$, and $f_2=
f-\e g$ in Example \ref{af}. Expand the formula for $\ti v_2$
given by \eqref{ah} in $\e$ to see
$$\ti v_2 = C_1\bpm 1\cr f(w)\epm + C_2\bpm
\overline{f(w)}\cr -1\epm   + \ O(\e),$$ where $w=x+\a
u+\a^{-1}v$, and
$$C_1=(1+|f(w)|^2), \quad C_2= (\bar{\a}-\a)((u-\a^{-2} v)f'(w)
+ g(w)).$$ Let $\hat \pi_2$ denote the projection onto the complex
line spanned by
$$\hat v_2= C_1\bpm 1\cr f(w)\epm + C_2\bpm
\overline{f(w)}\cr -1\epm.$$ Then as $\e\to 0$, $\ti\pi_2$ tends
to $\hat \pi_2$, the extended solution $\ti \psi$ tends to
$$\hat \psi= g_{\a, \hat \pi_2} g_{\a, \pi_1},$$ and the Ward
$2$-soliton tends to
$$\hat J= \bar b(\pi_1 + b\pi_1^\perp) (\hat \pi_2
+ b\hat \pi_2^\perp),
$$ where $b=\a/\bar \a$. In particular, if $\a=
i$, then we get an extended solution
\begin{equation}\label{am}
\hat\psi_2= g_{i, \hat\pi_2}g_{i, \pi_1},
\end{equation}
where $\pi_1(x,y,t)$ is the projection onto $\C\bpm 1\cr f(w)
\epm$ and $\hat \pi_2$ the projection onto the complex line
spanned by
$$\hat v_2= (1+|f(w)|^2)\bpm 1\cr f(w)\epm -2i (tf'(w)+ g(w))
\bpm \overline{f(w)}\cr -1\epm,$$ and $w= x+ iy$.  The limiting
Ward map $\hat J= -(\pi_1- \pi_1^\perp)(\hat \pi_2
-\hat\pi_2^\perp)$ is a $2$-soliton with non-trivial scattering,
and the extended solution $\hat\psi$ has a double pole at $\l=i$.
\eeg

Using similar limiting method, Ioannidou constructed extended Ward
$3$-solitons with a triple pole at $\l=i$ (see \cite{I1996}).

Below we apply algebraic BTs and an order $k$ limiting method to
construct $k$-solitons, whose extended solutions have pole data
$(z,k)$ for any $z\in \C\setminus\R$ and $k\geq 2$. To present
this method more clearly, we work on the $SU(2)$ case first. At
the end of this section, we will briefly explain how to generalize
this method to the $SU(n)$ case.

Let $z\in \C\setminus\R$ be a constant, and $\{
a_j(w)\}^\infty_{j=0}$ a sequence of rational functions in one
complex variable. Assume that $a_0(w)$ is not a constant function.
Let $a_j^{(i)}(w)$ denote the $i$-th derivative of $a_j$ with
respect to $w$. For any $\e\in\C$ with $|\e|$ small, let
\begin{align} \label{ww} w &=x+zu+z^{-1}v, \\ \notag
w_\e &=x+(z+\e)u +(z+\e)^{-1}v, \end{align} and $$f_{k,\e}
=\sum^k_{j=0} a_j(w_\e) \e^j.$$ A direct computation gives the
following Taylor expansions in $\e$:
\begin{enumerate}
    \item \begin{align*} w_\e -w &=\e u+((z+\e)^{-1} -z^{-1})v\\
    & =(u-z^{-2}v)\e +\sum_{l=2}^\infty z^{-l-1}v(-\e)^l.
    \end{align*}
    \item $$a_j(w_\e)=\sum_{l=0}^\infty
    \frac{a_j^{(l)}(w)}{l!}(w_\e-w)^l :=\sum_{l=0}^\infty
    b_{j,l}\e^l,$$ where  $b_{j,l}=b_{j,l}(x,u,v)$
    can be computed directly:
    \begin{equation*} \bca b_{j,0} &= a_j(w), \\
    b_{j,1} &= (u-z^{-2}v)a^\prime_j(w) \\ b_{j,2} &=
    \frac{(u-z^{-2}v)^2}{2} a_j^{\prime\prime}(w)
    +(z^{-3}v) a^\prime_j(w), \\ b_{j,3} &=
    \frac{(u-z^{-2}v)^3}{3!} a_j^{\prime\prime\prime}(w)
    +(u-z^{-2}v)(z^{-3}v) a_j^{\prime\prime}(w) -(z^{-4}v)a^\prime_j(w), \\
    b_{j,4} &= \frac{(u-z^{-2}v)^4}{4!} a^{(4)}_j (w)
    +\frac{a^{\prime\prime\prime}_j (w)}{2} (u-z^{-2}v)^2(z^{-3}v) \\
    & \quad +\frac{a^{\prime\prime}_j(w)}{2} (z^{-6}v^2
    -2(u-z^{-2}v)(z^{-4}v)) +(z^{-5}v) a^\prime_j (w), \\
    \cdots \eca \end{equation*}
    \item \begin{align*} f_{k,\e} &= \sum_{j=0}^k a_j (w_\e)\e^j
    =\sum_{j=0}^k \sum_{l=0}^\infty b_{j,l}\e^{j+l} \\
    &:= c_0+c_1\e+\cdots +c_k\e^k+O(\e^{k+1}), \end{align*} where
    $c_l =c_l(x,u,v) =\sum_{j=0}^l b_{j,l-j}$ are given below:
    \begin{align} \label{ao}
    \bca c_0 = a_0(w), \\
    c_1 = (u-z^{-2}v)a^\prime_0(w) +a_1(w), \\
    c_2 = \frac{(u-z^{-2}v)^2}{2}a^{\prime\prime}_0 (w)
    +(z^{-3}v)a^\prime_0(w) +(u-z^{-2}v)a^\prime_1(w) +a_2(w),\\
    c_3 = \frac{(u-z^{-2}v)^3}{3!} a^{\prime\prime\prime}_0 (w)
    +(u-z^{-2}v)(z^{-3}v)a^{\prime\prime}_0 (w) -(z^{-4}v)a^\prime_0
    (w) \\ \qquad+\frac{(u-z^{-2}v)^2}{2}a^{\prime\prime}_1 (w)
    +(z^{-3}v)a^\prime_1(w) +(u-z^{-2}v)a^\prime_2(w) +a_3(w),\\
    c_4 =\frac{(u-z^{-2}v)^4}{4!} a^{(4)}_0 (w)
    +\frac{a^{\prime\prime\prime}_0 (w)}{2} (u-z^{-2}v)^2(z^{-3}v) \\
    \qquad +\frac{a^{\prime\prime}_0(w)}{2} (z^{-6}v^2
    -2(u-z^{-2}v)(z^{-4}v)) +(z^{-5}v) a^\prime_0 (w) \\ \qquad
    +\frac{(u-z^{-2}v)^3}{3!} a^{\prime\prime\prime}_1 (w)
    +(u-z^{-2}v)(z^{-3}v)a^{\prime\prime}_1 (w) -(z^{-4}v)a^\prime_1
    (w) \\ \qquad+\frac{(u-z^{-2}v)^2}{2}a^{\prime\prime}_2 (w)
    +(z^{-3}v)a^\prime_2(w) +(u-z^{-2}v)a^\prime_3(w) +a_4(w), \\
    \cdots \eca \end{align}
\end{enumerate}

From the above computation, we see that $c_j$'s are rational
functions in $x, u, v$, hence are rational in $x,y$ and $t$.
Singularities of $c_j$ consist of finitely many straight lines in
$\R^{2,1}$ given by $w=x+zu+z^{-1}v=p_1,p_2,\cdots$, where
$p_1,p_2,\cdots$ are the poles of the rational functions
$a_0,\cdots,a_j$.

For $k\geq 1$, let
$$v_{k,\e}= \bpm 1\\ f_{k-1,\e} \epm,$$
$\pi_{k,\e}$ the Hermitian projection of $\C^2$ onto $\C
v_{k,\e}$. Define $\psi_{k,\e}$ and $\hat \psi_k$ by induction as
follows:
\begin{align*}
&\psi_{1,\e} =g_{z+\e,\pi_{1,\e}}, \quad \text{and}
    \quad \hat\psi_1= \lim_{\e\to 0} \psi_{1,\e},\\
& \psi_{k,\e}= g_{z+\e, \pi_{k,\e}}\ast \hat\psi_{k-1},
    \quad \text{and} \quad \hat\psi_{k} =\lim_{\e\to 0} \psi_{k,\e}.
\end{align*}
 Let
$$\ti v_{k,\e}= \hat \psi_{k-1}(z+\e)(v_{k,\e}),$$
By Theorem \ref{thm:BT},
$$\psi_{k,\e}= g_{z+\e,\ti\pi_{k,\e}}\hat\psi_{k-1},$$
where $\ti\pi_{k,\e}$ is the projection onto $\C \ti v_{k,\e}$.

\begin{thm} \label{thm:LP}
Let $a_0, a_1, \ldots $ be a sequence of rational functions from
$\C$ to $\C$, and let $v_{k,\e}, \pi_{k,\e}, \psi_{k,\e},
\hat\psi_k$ and $\ti v_{k,\e}$ be defined as above.  Then we have
\begin{enumerate}
    \item
$\ti v_{k,\e} =\hat v_k +\e y_{k,1}
    +\e^2 y_{k,2} +\cdots$, where
\begin{subequations}
\begin{gather}
 \hat v_k = \bpm 1\\c_0 \epm
    +\sum_{j=1}^{k-1} (z-\bar z)^j P_{k-1,j} \bpm 0\\c_j \epm, \label{q1} \\
  P_{l , j} =\sum_{l\ge i_1  >\cdots  >i_j\ge 1} \hat \pi_{i_1}^\perp \cdots
    \hat \pi_{i_j}^\perp.   \label{kk}
 \end{gather}
\end{subequations}
 Moreover, all entries of $\hat v_k$ are rational functions
in $x,y$ and $t$.
 \item
 $\hat\psi_k= g_{z, \hat\pi_k}\cdots g_{z, \hat\pi_1}$
    is a minimal factorization and is an extended Ward map with only
a pole at $\l= z$ of multiplicity
    $k$, where $\hat \pi_k$ is the Hermitian projection of $\C^2$ onto $\C
    \hat v_k$. Moreover, $\hat \pi_k$ is smooth, and for each  fixed
$t$, $\lim_{\n x\n^2 +\n y\n^2 \to \infty} \hat \pi_k(x,y,t)$
    exists.
 \item
 The Ward map associated to $\hat \psi_k$,
    $$J_k = \frac{1}{|z|^k}(\bar z\hat\pi_1 +z\hat\pi_1^\perp)
    (\bar z\hat\pi_2 +z\hat\pi_2^\perp) \cdots (\bar z\hat\pi_k
    +z\hat\pi_k^\perp),$$ is smooth and satisfies the boundary
    condition \eqref{bc1}, and all entries of $J_k$ are rational
    functions in $x,y$ and  $t$.
\end{enumerate}
\end{thm}

\begin{proof} We prove the Theorem by induction on $k$. For $k=1$,
Theorem is clearly true. Suppose the Theorem is true for $k$. We
will prove that (1)--(3) hold for $k+1$.

 \ms
\ni (1) By Theorem \ref{thm:BT} and induction hypothesis, we have
\begin{align*}
&\ti v_{k+1,\e} = \hat \psi_k(z+\e) v_{k+1,\e}
\\ &=({\rm I} +\frac{z-\bar z}{\e} \hat{\pi}^\perp_k)\hat\psi_{k-1}(z+\e)
\begin{pmatrix} 1 \\ f_{k,\e} \end{pmatrix} \\
&= ({\rm I} +\frac{z-\bar z}{\e}
\hat{\pi}^\perp_k)\hat\psi_{k-1}(z+\e)
\begin{pmatrix} 1 \\ f_{k-1,\e} +a_k(w_\e) \e^k
\end{pmatrix} \\ &= ({\rm I} +\frac{z-\bar z}{\e}
\hat{\pi}^\perp_k) \left( \hat\psi_{k-1}(z+\e)
\begin{pmatrix} 1 \\ f_{k-1,\e} \end{pmatrix} +\hat\psi_{k-1}(z+\e)
\begin{pmatrix} 0 \\ a_k(w_\e) \e^k \end{pmatrix} \right) \\
&= ({\rm I} +\frac{z-\bar z}{\e} \hat{\pi}^\perp_k)\left( \hat v_k
+\e y_{k,1} +\e (z-\bar
z)^{k-1} P_{k-1,k-1} \begin{pmatrix} 0 \\
a_k(w) \end{pmatrix} +O(\e^2) \right)\\ &= \hat v_k + (z-\bar z)
\hat\pi_k^\perp \left( y_{k,1} +(z-\bar z)^{k-1} P_{k-1,k-1}
\begin{pmatrix} 0 \\ a_k(w) \end{pmatrix} \right) +O(\e).
\end{align*} In the last step we have used $\hat v_k \in \Im
\hat\pi_k$. Therefore all terms of negative powers of $\e$ vanish
in the Laurent series expansion of $\ti v_{k+1,\e}$ in $\e$.

  The Laurent series expansion of $\hat \psi_k(z+\e)$ in $\e$ is
  \begin{align} \notag \hat \psi_k(z+\e) &= (\I +\frac{z-\bar z}{\e}
\hat \pi_k^\perp) \cdots (\I +\frac{z-\bar z}{\e} \hat \pi_1^\perp) \\
\label{as} &= \I +\frac{z-\bar z}{\e} P_{k,1} +\cdots +
\frac{(z-\bar z)^k}{\e^k} P_{k,k}.
  \end{align}
Substituting \eqref{as} to $\ti v_{k+1,\e}$ and using the fact
that the Laurent series expansion of $\ti v_{k+1,\e}$ has no
$\e^{-j}$ terms with $j>0$,  we have
\begin{align*} \ti v_{k+1,\e} &= \hat \psi_k(z+\e) v_{k+1,\e}
=\hat \psi_k(z+\e) \bpm 1\\ f_{k,\e} \epm \\
&= \hat \psi_k(z+\e) \left( \bpm 1\\c_0\epm +\e \bpm 0\\c_1
\epm +\cdots +\e^k \bpm 0\\c_k \epm +O(\e^{k+1}) \right) \\
&=\bpm 1\\c_0\epm +(z-\bar z)P_{k,1} \bpm 0\\c_1 \epm +\cdots
+(z-\bar z)^k P_{k,k} \bpm 0\\c_k \epm +O(\e).
  \end{align*}
Therefore we obtain
  \begin{equation}\label{bu}
  \hat v_{k+1}= \bpm 1\\ c_0\epm + \sum_{j=1}^k (z-\bar z)^j P_{k,
j}\bpm 0\\ c_j\epm.
  \end{equation}

  By induction hypothesis, $\hat
\pi_1,\cdots, \hat \pi_k$ are smooth, and all of their entries are
rational in $x,y$ and $t$. Thus $P_{k,1},\cdots, P_{k,k}$ have the
same analytic properties as $\pi_j$'s. Together with the analytic
properties of $c_j$'s, we see that all entries of $\hat v_{k+1}$
are rational in $x,y$ and $t$.

  \ms

\ni (2) By (1), we have $$\lim_{\e\to 0} \ti \pi_{k+1,\e} =\hat
\pi_{k+1},$$ where $\hat \pi_{k+1}$ is the projection onto $\C
\hat v_{k+1}$. Since all entries of $\hat v_{k+1}$ are rational in
$x,y$ and $t$, $\hat \pi_{k+1}$ is smooth and for each fixed $t$,
$\lim_{|(x,y)| \to \infty} \hat \pi_{k+1}(x,y,t)$
exists.

  Next we claim that $\hat \psi_{k+1}$ is an extended Ward map with a
pole at $\l=z$ of multiplicity $k+1$.  To see this, first note
that
  \begin{align*}
  \hat \psi_{k+1} &=\lim_{\e\to 0}
g_{z+\e,\pi_{k+1,\e}} \ast \hat \psi_k = \lim_{\e\to 0}
g_{z+\e,\ti \pi_{k+1,\e}} \hat \psi_k \\ &= g_{z,\hat \pi_{k+1}}
g_{z,\hat \pi_k}\cdots g_{z,\hat \pi_1}= (\I+\frac{z-\bar
z}{\l-z}\hat \pi_{k+1})\cdots (\I+\frac{z-\bar z}{\l-z} \hat\pi_1) \\
  &= \I + \sum_{j=1}^{k+1} \frac{(z-\bar z)^j}{(\l-z)^j} \ P_{k+1,j}.
  \end{align*}
  By Theorem~\ref{thm:BT}, $g_{z+\e,\pi_{k+1,\e}} \ast \hat \psi_k$ is
an extended Ward map for small $|\e |  >0$.  By continuity, so is
$\hat \psi_{k+1}$.  The coefficient of $(\l-z)^{-k-1}$ of $\hat
\psi_{k+1}$ is $(z-\bar z)^{k+1} P_{k+1, k+1}$.
  To show that $\hat \psi_{k+1}$ has a pole
at $\l =z$ of multiplicity $k+1$, it suffices to show that
$P_{k+1,k+1} =\hat \pi_{k+1}^\perp \hat \pi_k^\perp \cdots \hat
\pi_1^\perp \not= 0$. For this purpose, we write

  \begin{equation}\label{dn}
  P_{k,j}=P_{k-1,j} +\hat \pi_k^\perp P_{k-1,j-1}.
   \end{equation}
  So \eqref{bu}
for $k+1$ can be written as
  \begin{align} \hat v_{k+1} &=
\bpm 1\\c_0 \epm +\sum_{j=1}^{k-1} (z-\bar z)^j P_{k-1,j} \bpm
0\\c_j \epm + \hat \pi_k^\perp \sum_{j=0}^{k-1} (z-\bar z)^{j+1}
P_{k-1,j} \bpm 0\\c_{j+1} \epm \\ \label{do} &= \hat v_k +(z-\bar
z) \hat \pi_k^\perp \sum_{j=0}^{k-1} (z-\bar z)^j P_{k-1,j} \bpm
0\\c_{j+1} \epm.
\end{align}
  By the induction hypothesis, $\hat \psi_k= g_{z,\hat\pi_k}\cdots
g_{z,\hat\pi_1}$ is a minimal factorization.  So
$(\Im\pi_{j+1})\cap (\Im\pi_j^\perp)=0$ for $1\le j\le k-1$.  By
induction hypothesis $\hat v_k\not=0$.  Formula \eqref{do} implies
that $\Im\hat\pi_{k+1}\cap \Im\hat \pi_k^\perp=0$.  Hence $\hat
\psi_{k+1}$ is a minimal factorization.  By Proposition~\ref{cr},
$\hat \pi_{k+1}^\perp \hat \pi_k^\perp \cdots \hat \pi_1^\perp
\not= 0$. Thus (2) holds for $k+1$.

\ni (3) The expression for $J_{k+1} =\hat \psi_{k+1}(0)^{-1}$ is
straightforward. Since all $\pi_j$'s are smooth, and for each
fixed $t$, $\lim_{\n x\n^2 +\n y\n^2 \to \infty} \pi_j(x,y,t)$
exist, $J_{k+1}$ is also smooth and satisfies the boundary
condition \eqref{bc1}. The entries of $J_{k+1}$ are rational in
$x,y$ and $t$ because all $\pi_j$'s have this property.
\end{proof}

We give some explicit formulas for $\hat v_k$ with $k$ small: \ben
  \item $\hat v_1 =\begin{pmatrix} 1\\ c_0 \end{pmatrix}$.
  \item  $\hat v_2 =\hat v_1 +(z-\bar z) \hat\pi_1^\perp
  \begin{pmatrix} 0\\ c_1 \end{pmatrix}$.
  The corresponding Ward map coincides with the one in
  Example~\ref{ag} if $z=\a$ and $(a_0,a_1)=(f,g)$.
  \item $\hat v_3 = \hat v_2 +(z-\bar z)\hat\pi_2^\perp
  \left( \bpm 0\\c_1 \epm +(z-\bar z) \hat\pi_1^\perp
  \begin{pmatrix} 0\\ c_2 \end{pmatrix} \right)$.
  \item
\begin{align*} & \hat v_4- \hat v_3 \\
& \ = (z-\bar z)\hat\pi_3^\perp \left(
 \bpm 0\\ c_1 \epm +(z-\bar z) (\hat\pi_2^\perp +\hat\pi_1^\perp)
 \bpm 0\\c_2 \epm
  +(z-\bar z)^2 \hat\pi_2^\perp \hat\pi_1^\perp \bpm 0\\c_3 \epm
  \right).
\end{align*}
\een

\ms

  We briefly explain how to construct Ward solitons with pole data
$(z,k)$ for the $SU(n)$ case next. Choose a sequence of rational
maps $a_j:\C \to \C^n$, $j=0,1,2,\cdots$. Let $\pi_{j,\e}$ be the
Hermitian projection of $\C^n$ onto $\C \sum_{l=0}^{j-1} a_l(w_\e)
\e^l$, $w_\e=x+(z+\e)u+(z+\e)^{-1}v$, for $1\le j\le k$. Then the
same computation and proof as in the $SU(2)$ case imply that $\hat
\psi_k =g_{z,\hat \pi_k} g_{z,\hat \pi_{k-1}}\cdots g_{z,\hat
\pi_1}$ is a minimal factorization and is an extended Ward map
with pole data $(z,k)$.

  Note that
all $\hat \pi_j$'s are of rank one in the above construction.  But
the same limiting method also produces extended Ward solitons of
the form
  $$g_{z,\pi_k}\cdots g_{z,\pi_1}$$
with $\rk(\pi_1)\geq \cdots\geq \rk(\pi_k)$.  To see this, let
$n-1\geq n_1\geq \cdots \geq n_k\geq 1$ be integers, and $a_{i1},
\ldots, a_{i n_i}$ $\C^n$-valued rational maps on $\C$ for $1\leq
i\leq k$. Suppose
  $$a_{11}\wedge\cdots \wedge a_{1n_1} \not= 0$$ generically, and  let
$\pi_1(w)$ denote the projection of $\C^n$ onto the linear span of
$a_{11}(w), \ldots, a_{1n_1}(w)$.   Let $w=x+zu+z^{-1}v$, and
$w_\e= x+ (z+\e)u +(z+\e)^{-1}v$.   Then
  \begin{align*}
  \hat v_{2i}&:=\lim_{\e\to 0} g_{z,\pi_1}(z+\e)(a_{1i}(w_\e)+ \e
a_{2i}(w_\e))\\
  &= a_{1i}(w)+ (z-\bar z)\pi_1^\perp((u-z^{-2}v) a_{1i}'
(w)+a_{2i}(w)), \quad 1\leq i\leq n_2.
  \end{align*}
  Let $\pi_2$ denote the projection onto the linear span of
$\hat v_{21}(w), \ldots, \hat v_{2n_2}(w)$.  Then $g_{z,
\pi_2}g_{z,\pi_1}$ is an extended solution, and $\rk(\pi_2)= n_2$.
It is easy to see that $\Im\pi_2\cap \Im\pi_1^\perp=0$. Hence
$g_{z,\pi_2}g_{z,\pi_1}$ is a minimal factorization. Similar
computations give the construction of extended Ward maps with pole
data $(z, k)$ and arbitrary rank data $(n_1, \ldots, n_k)$.

\bs

\section{Ward solitons with general pole data}

   We associate to each extended Ward map with pole data $(z,k)$ a
generalized algebraic BT.  Use these generalized BTs, we construct
extended Ward maps that have general pole data $(z_1, \ldots, z_k,
n_1,\ldots, n_k)$.

We first give a more general algebraic BT (Theorem \ref{thm:BT}):

 \bthm  \label{thm:GBT}
Let $\phi$ be an extended Ward map with pole data $(z,k)$, and
$\psi$ an extended Ward map that is holomorphic and non-degenerate
at $\l =z,\bar z$.  Then there exist unique $\ti \phi$ and
$\ti\psi$ such that $\ti \phi \psi =\ti \psi \phi$, where
$\ti\phi$ has pole data $(z,k)$, and $\ti\psi$ is holomorphic and
non-degenerate at $\l=z,\bar z$.  Moreover,
\begin{equation}
\psi_k= \ti \phi \psi =\ti \psi \phi,
\end{equation}
is a new extended Ward map and $\ti\phi, \ti \psi$ are constructed
algebraically. \ethm

\begin{proof}
 It follows from Theorem \ref{bf} that we can factor $\phi$ as
product of $k$ simple elements.
 $$\phi=g_{z,\pi_k}\cdots g_{z,\pi_1}.$$
  Let $\ti \pi_1$ be the projection onto
$\psi(z)\Im \pi_1$, and $\ti \psi_1=g_{z,\ti \pi_1} \psi
g_{z,\pi_1}^{-1}$. Residue calculus implies that $\ti \psi_1$ is
holomorphic at $\l=z,\bar z$. For $j=2,\cdots,k$, we
 define $\ti \pi_j$ and $\ti \psi_j$ recursively by
$$\Im \ti \pi_j =\ti \psi_{j-1}(z)\, \Im\pi_j, \quad \text{and}
\quad \ti \psi_j =g_{z,\ti \pi_j}\ti \psi_{j-1}
g_{z,\pi_j}^{-1}.$$ Again $\ti \psi_j$ is holomorphic at
$\l=z,\bar z$ for $j=2,\cdots,k$ by residual calculus. Let $\ti
\phi=g_{z,\ti \pi_k}\cdots g_{z,\ti \pi_1}$, and $\ti \psi =\ti
\psi_k$.  By construction,  $\ti \phi \psi =\ti \psi \phi$. Next
we prove uniqueness.  Suppose $\hat\phi$ has pole data $(z,k)$ and
$\hat \psi$ is holomorphic and non-degenerate at $z, \bar z$ and
$\hat \phi\psi= \hat \psi\phi$. Then $\phi\psi^{-1}= \ti
\psi^{-1}\ti\phi = \hat \psi^{-1}\hat \phi$. So we have
$$\hat \psi\ti\psi^{-1}= \hat \phi \ti\phi^{-1}.$$
But the left hand side is holomorphic at $z, \bar z$ and the right
hand side is holomorphic at $\l\in \C\setminus \{z, \bar z\}$ and
is equal to $\I$ at $\l=\infty$.  Hence it must be the constant
identity.  This proves $\ti\phi=\hat \phi$ and $\ti \psi=\hat
\psi$.

The same proof of Theorem \ref{thm:BT} implies that $\psi_k=
\ti\psi\phi= \ti\phi\psi$ is an extended Ward map.
\end{proof}

We use $\phi\ast \psi$ to denote the new extended solution
$\psi_1= \ti\phi \psi$ constructed in the above Theorem, and call
$$\psi\mapsto \phi\ast \psi$$ the {\it generalized B\"acklund
transformation generated by $\phi$\/}.

The proof of Theorem \ref{thm:GBT} implies that if $\psi$ and
$\phi$ are extended Ward maps with pole data $(z_1, n_1)$ and
$(z_2, n_2)$ respectively and $z_1\not= z_2, \bar z_2$, then
$\phi\ast \psi= \psi\ast \phi$. Same argument gives the following
Corollary:

\bcor Let $z_1, \ldots, z_r\in \C\setminus \R$ such that $z_i\not=
z_j, \bar z_j$ for all $i\not=j$, and $\phi_j$ an extended Ward
map with pole data $(z_j, n_j)$ for $1\leq j\leq r$.  Let $\sigma$
be a permutation of $\{1, \ldots, k\}$.  Then
$$\phi_{\sigma(1)}\ast (\phi_{\sigma(2)}\ast ( \cdots
\ast\phi_{\sigma(r)})\cdots ) = \phi_1\ast (\phi_2 \ast (\cdots
\ast \phi_r)\cdots).
 $$
\ecor

\beg \label{bk} Extended Ward $4$-solitons into $SU(2)$ with two
double poles.
\par

Choose $z_1,z_2\in\C \setminus \R$ with $z_1\not=z_2, \bar z_2$,
and rational functions $a_0(w),a_1(w)$ and $b_0(w),b_1(w)$. Let
$w_i=x+z_iu +z_i^{-1}v$, $i=1,2$. By the construction of section
5, we have two extended Ward solitons
$$\phi=g_{z_2,\pi_2}g_{z_2,\pi_1},\quad \psi=g_{z_1,\tau_2}
g_{z_1,\tau_1},$$ where $$\Im \pi_1=\C v_1 =\C \bpm 1\\b_0(w_2)
\epm,\quad \Im \pi_2=\C v_2,$$  \begin{align*} v_2 &=\bpm
1\\b_0(w_2) \epm +(z_2-\bar z_2)\pi_1^\perp \bpm 0\\
(u-z_2^{-2}v)b_0^\prime(w_2) +b_1(w_2) \epm \\   &=\bpm
1\\b_0(w_2) \epm -\frac{(z_2-\bar
z_2)((u-z_2^{-2}v)b_0^\prime(w_2) +b_1(w_2))}{1+|b_0(w_2)|^2} \bpm
\overline{b_0(w_2)}\\ -1 \epm,
\end{align*} and $$\Im \tau_1=\C q_1 =\C \bpm 1\\a_0(w_1)
\epm,\quad \Im \tau_2=\C q_2,$$ \begin{align*} q_2 &=\bpm
1\\a_0(w_1) \epm +(z_1-\bar z_1) \tau_1^\perp \bpm 0\\
(u-z_1^{-2}v)a_0^\prime(w_1) +a_1(w_1) \epm \\  &=\bpm 1\\a_0(w_1)
\epm -\frac{(z_1-\bar
z_1)((u-z_1^{-2}v)a_0^\prime(w_1)+a_1(w_1))}{1+|a_0(w_1)|^2} \bpm
\overline{a_0(w_1)}\\ -1 \epm.\end{align*}  Apply generalized
B\"acklund transformation $\phi\ast \psi$ to get an extended
solution with two double poles at $z_1, z_2$. By
Theorem~\ref{thm:GBT}, $$\phi \ast\psi =g_{z_2,\ti \pi_2}
g_{z_2,\ti \pi_1} g_{z_1,\tau_2} g_{z_1,\tau_1},$$ where $\Im \ti
\pi_1 =\C\ti v_1=\C \psi(z_2) v_1$, and $\Im \ti \pi_2 =\C \ti v_2
=\C \ti \psi_1 (z_2) v_2$. Compute the following limit
$$(z_2-\bar z_2) \ti \psi_1 (z_2) =\lim_{\l\to z_2} (\l-\bar z_2)
g_{z_2,\ti \pi_1}(\l) \psi(\l) g_{z_2,\pi_1} (\l)^{-1}$$ to see
$$\ti \psi_1 (z_2) =\psi(z_2)\pi_1 +\ti \pi_1^\perp \left(
\psi(z_2) +(z_2-\bar z_2) \frac{\partial\psi (z_2)}{\partial \l}\,
\pi_1\right).$$ The associated Ward $4$-soliton is given by
$$J_4= \frac{1}{|z_1z_2|^2} (\bar z_1 \tau_1 +z_1 \tau_1^\perp)
(\bar z_1 \tau_2 +z_1 \tau_2^\perp) (\bar z_2 \ti \pi_1 +z_2 \ti
\pi_1^\perp) (\bar z_2 \ti \pi_2 +z_2 \ti \pi_2^\perp).$$ \eeg

\ms

Let
$$\C_\pm=\{a\pm ib\n b>0\}$$
denote the upper and lower half plane of $\C$. We claim that to
construct Ward solitons with general pole data, we may assume all
the poles lie in the upper half plane $\C_+$.  This claim follows
from two remarks below:

(1) A direct computation implies that
$$\frac{\l-z}{\l-\bar z}\ g_{z,\pi} (\l) = \I + \frac{\bar
z-z}{\l-\bar z}\pi = g_{\bar z, \pi^\perp}.$$

(2) Let $\psi$ be an extended Ward soliton.  By Theorem \ref{bf},
we can factor
$$\psi= g_{z_1, \pi_1}\cdots g_{z_r, \pi_r}.$$
Suppose $z_1, \ldots, z_k\in \C_-$ and the rest of the poles lie
in $\C_+$.  Let
$$f(\l)= \prod_{j=1}^k \frac{\l- z_{j}}{\l-\bar z_{j}}.$$
Then $f\psi$ still satisfies \eqref{eq:auxiliary}, hence is an
extended solution.  But $f\psi$ has poles at $\bar z_1, \ldots,
\bar z_k, z_{k+1}, \ldots, z_r$, which all lie in $\C_+$.  The
Ward maps corresponding to $\psi$ and $f\psi$ are $J=\psi(0)^{-1}$
and $J_1= \psi(0)^{-1}/f(0)$ respectively.  Note $f(0)$ is a
constant complex number of length $1$.  So we do not lose any Ward
maps by assuming that all poles lie in $\C_+$.

\bcor Given distinct  $z_1,\cdots,z_r \in\C_+$ and positive
integers $n_1$, $\ldots$, $n_r$,  there is a family of Ward
solitons whose extended solutions have pole data
$(z_1,\cdots,z_r,n_1,\cdots,n_r)$. \ecor

\begin{proof}  Let
$\phi_{z_j,n_j} =g_{z_j,\pi_{n_j}}\cdots g_{z_j,\pi_1}$ be an
extended Ward soliton with pole data $(z_j,n_j)$ constructed in
section 5. Apply Theorem \ref{thm:GBT} repeatedly to $\phi_{z_1,
n_1}$ to get the extended solution
$$\phi= \phi_{z_r,n_r}\ast (\cdots\ast(\phi_{z_2,n_2}\ast
\phi_{z_1,n_1})\cdots ).$$ Then $\phi$  has pole data
$(z_1,\cdots,z_r,n_1,\cdots,n_r)$.
\end{proof}

In the rest of the section, we prove that a general Ward soliton
can be constructed by applying generalized B\"acklund
transformations to an extended Ward soliton with pole data
$(z,k)$.

\bthm \label{bh} Suppose $\psi$ is an extended solution of the
Ward equation, and $\psi= f_1f_2$, such that \ben \item $\l\mapsto
f_2(x,y,t,\l)$ is an element of the group $\cs_r(S^2,GL(n))$ and
has poles only at $z_1, \ldots, z_k$, \item  $f_1$ is holomorphic
and non-degenerate at $\l=z_1, \ldots, z_k$ and $\bar z_1, \ldots,
\bar z_k$. \een Then $f_2$ is also an extended solution of the
Ward equation. \ethm

\begin{proof}
Let $Pf=\l f_x- f_u$, $Qf= \l f_v - f_x$, and $D=\{z_1, \ldots,
z_k, \bar z_1, \ldots, \bar z_r\}$.     Use $f_2= f_1^{-1}\psi$ to
compute directly to get
$$A_2:=(Pf_2)f_2^{-1} = P(f_1^{-1})f_1 + f_1^{-1}(P\psi)\psi^{-1} f_1.$$
Since $\psi$ is an extended solution, $A=(P\psi)\psi^{-1}$ is
independent of $\l$. Because $f_1$ is holomorphic and
non-degenerate at points in $D$, the right hand side of $A_2$ is
holomorphic at points in $D$.   But $f_2$ is assumed to be
holomorphic in $\C\setminus D$, so $(Pf_2)f_2^{-1}$ is holomorphic
for $\l\in\C\setminus D$.  So $A_2(x,y,t,\l)$ is holomorphic for
all $\l\in \C$.  But $f_2\in \cs_r(S^2,GL(n))$ implies that
$(Pf_2)f_2^{-1}$ is holomorphic at $\l=\infty$.  Hence $A_2$ is
independent of $\l$.  Similarly, $(Qf_2)f_2^{-1}$ is independent
of $\l$.  This proves that $f_2$ is an extended solution of the
Ward equation.
\end{proof}

\bcor \label{bo} Suppose $\psi= \phi_1\cdots \phi_r$ is an
extended Ward soliton such that $\phi_j$ has pole data $(z_j,n_j)$
and $z_1,\ldots, z_r\in\C_+$ are distinct.  Then for $1\le j\le
r-1$,  \ben \item $\psi_j= \phi_{j+1}\ldots \phi_r$ is also an
extended Ward soliton, \item there exists a unique extended Ward
soliton $\ti \phi_j$ with pole data $(z_j,n_j)$ so that
$\psi_{j-1}= \ti\phi_j\ast \psi_j$, \item $\psi$ can be
constructed by applying the generalized algebraic B\"acklund
transformations repeatedly to $\phi_r$. \een \ecor

\bprop \label{bl} If $\psi$ is an extended Ward soliton with pole
data $(z_1, \ldots, z_r,$ $ n_1, \ldots, n_r)$, then for $1\le
j\le r$, there exists a unique extended Ward soliton $\phi_j$ with
pole data $(z_j, n_j)$ so that
$$\psi= \phi_1 \ast(\phi_2\ast (\cdots \ast\phi_r)\cdots ).$$
\eprop

\begin{proof}
By Uhlenbeck's factorization Theorem \ref{bf} and the
permutability Theorem \ref{bg}, we can factor $\psi$ as
\begin{equation} \label{bi}
\psi=f_2\cdots f_r f_1 = g_1(g_2\cdots g_r)
\end{equation}
 such that $f_j, g_j\in \cs_r(S^2,GL(n))$ have pole data
$(z_j, n_j)$.  We prove the Proposition  by induction on $r$.  If
$r=1$, the Proposition is automatically true. Suppose the
Proposition is true for $r=n-1$.  Then by Theorem \ref{bh} both
$f_1$ and $g_2\cdots g_r$ are extended solutions.  By induction
hypothesis, there exist extended Ward maps $h_2, \ldots, h_r$ with
pole data $(z_2, n_2), \ldots, (z_r, n_r)$ respectively such that
$$g_2\cdots g_r= h_2\ast (h_3\ast (\cdots \ast h_r)\cdots).$$
Equation \eqref{bi} implies that $\psi = f_1\ast (g_2\cdots g_r)$,
which is equal to
$$f_1\ast (h_2\ast (\cdots \ast h_r)\cdots ).$$
\end{proof}

\bs

\section{Analytic BT and Holomorphic vector bundles}

  In this section, we  generalize some of Uhlenbeck's results on
unitons to Ward solitons.   In particular, we
  \ben
   \item  derive the analytic B\"{a}cklund transformation
BT$_{z,\psi}$ (given in the introduction) for the Ward equation,
  \item  associate to each Ward soliton and complex number $z\in
\C\setminus \R$ a holomorphic structure on the trivial
$\C^n$-bundle over $S^2$,
  \item  prove that to find solutions $\pi$ of  BT$_{z,\psi}$ is
equivalent to find a one parameter family of holomorphic
subbundles of the  $\C^n$-bundle with respect to the holomorphic
structure given in (2) that satisfy certain first order PDE
system.
  \een

Suppose $\psi(\lambda)(x,y,t)=\psi(x,y,t,\l)$ is an extended
solution of the Ward equation with $A=(\l\psi_x-\psi_u)\psi^{-1}$,
and $B=(\l\psi_v-\psi_x)\psi^{-1}$. Motivated by the construction
of B\"acklund transformations for soliton equations, we seek a new
extended solution of the form $\psi_1(\lambda) =g_{z,\pi}(\lambda)
\psi(\lambda)$ for some smooth map $\pi$ from $\R^{2,1}$ to the
space of rank $k$ Hermitian projections of $\C^n$.  The condition
$\psi_1$ satisfies \eqref{eq:auxiliary} implies that

\begin{equation}\label{dh}
\bca (\l g_x - g_u)g^{-1} + gAg^{-1}= \ti A , &\\
(\l g_v -g_x)g^{-1} +g Bg^{-1}=\ti B ,\eca
\end{equation}
 for some $\ti A$ and $\ti B$ independent of $\l$,
where $g=g_{z, \pi}= \I+\frac{z-\bar z}{\l-z}\pi^\perp$. The
reality condition implies that $g_{z,\pi}(\l)^{-1}= \I +\frac{\bar
z-z}{\l-\bar z}\pi^\perp$. Since the right hand side of \eqref{dh}
is holomorphic in $\l\in \C$, the residue of the left hand side at
$\l= z$ must be zero, which gives
$$\bca
( z \pi^\perp_x-\pi^\perp_u)\pi + \pi^\perp A\pi=0,\\
( z\pi_v^\perp-\pi_x^\perp)\pi +\pi^\perp B\pi=0. \eca$$ But
$\pi^\perp\pi=0$ implies that
$$-\pi^\perp(z\pi_x-\pi_u) = ( z \pi^\perp_x-\pi^\perp_u)\pi.$$
So we get the following system of first order partial differential
equations for $\pi$:
\begin{equation} \left\{ \begin{array}{l}
\pi^\perp (z\pi_x -\pi_u - A\pi) =0, \\
\pi^\perp (z\pi_v -\pi_x -B\pi) =0. \\
\end{array} \right. \label{singBT} \end{equation}
The residue at $\l= \bar z$ is zero gives the same system
\eqref{singBT}, which is the analytic B\"acklund transformation
(BT) for the Ward equation.   So we have proved the following:

\bprop \label{cy} Let $\psi$ be an extended solution of the Ward
equation, $A=(\l\psi_x-\psi_u)\psi^{-1}$, and
$B=(\l\psi_v-\psi_x)\psi^{-1}$. Given a smooth map $\pi$ from $
\R^{2,1}$ to the space of Hermitian projections, $g_{z,\pi}\psi$
is an extended solution of the Ward equation if and only if $\pi$
satisfies \eqref{singBT}. \eprop

 Ioannidou and Zakrzewski proved the above Proposition for $z=i$ in
\cite{IZ1998}.  However, no general solutions of the system
\eqref{singBT} were given.

 As a consequence of Theorem \ref{thm:BT} and Proposition \ref{cy} we have

\bprop Suppose $\psi(x,u,v,\l)$ is an extended Ward map and is
holomorphic and non-degenerate at $\l=z$, and $f_1, \ldots, f_k$
are meromorphic maps from $\C$ to $\C^n$ that are linearly
independent except at finitely many points of $\C$.  Let $w=
x+zu+z^{-1}v$, $\ti f_i= \psi(\cdots, z)(f_i)$, and $\ti\pi$ be
the Hermitian projection of $\C^n$ onto the span of $\ti f_1,
\ldots, \ti f_k$. Then $\ti\pi$ is a solution of ${\rm
BT\/}_{z,\psi}$, or equivalently, $g_{z,\ti\pi}\psi$ is again an
extended solution. \eprop

Thus if $\psi$ is holomorphic and non-degenerate at $\l= z$, then
the above Proposition gives an algebraic method to construct
solutions of BT$_{z,\psi}$.

When $z=i$, $\pi$ is independent of $t$, and $\psi$ is an extended
$k$-uniton, \eqref{singBT} is the singular BT used by Uhlenbeck in
\cite{U1989} to add one more uniton to the given $k$-uniton. She
also proved that a solution of \eqref{singBT} for a uniton can be
interpreted as a holomorphic subbundle  that satisfies an
algebraic constraint.
 In this section, we show that a solution $\pi$ of the analytic BT
\eqref{singBT} for Ward map can also be interpreted in terms of
holomorphic subbundle, but it now must satisfy a first order PDE
constraint. We explain this next.

If we make a suitable linear change of coordinates of $\R^{2,1}$,
then the operator $z\p_x-\p_u$ becomes a $\bar \p$ operator.  To
see this, let $z=\a+i\b$.  Then
\begin{equation*}
w= x+zu+z^{-1}v=(x+ \a u+\frac{\a v}{\a^2+\b^2}v)+ i(\b u-\frac{\b
v}{\a^2+\b^2}).
\end{equation*}
Make a coordinate change:
\begin{equation}\label{el}
\bca
p= x+ \a u +\frac{\a v}{\a^2+\b^2},&\\
q= \b u-\frac{\b v}{\a^2+\b^2},&\\
r=  v. \eca
\end{equation}
 Then
$$\bca
\p_x= \p_p,&\\
\p_u= \a\p_p+ \b\p_q,&\\
\p_v= \frac{\a}{\a^2+\b^2}\p_p-\frac{\b}{\a^2+\b^2}\p_q +\p_r.
\eca$$

A direct computation gives
\begin{align*}
z\p_x-\p_u&= (\a+i\b)\p_p -(\a\p_p+\b\p_q)= i\b\p_p -\b \p_q \\
&= i\b(\p_p+i\p_q) = 2i\b \p_{\bar w}.\\
z\p_v-\p_x&=
(\a+i\b)(\frac{\a}{\a^2+\b^2}\p_p-\frac{\b}{\a^2+\b^2}\p_q +\p_r)
-\p_p\\
&=\frac{2i\b}{\a-i\b}\left(\p_{\bar w} -
\frac{\a^2+\b^2}{2i\b}\p_r \right)
\end{align*}
Use $z=\a +i\b$ and substitute the above formulas into the
analytic BT \eqref{singBT} to get
\begin{equation}\label{cu}
\bca
\pi^\perp(\p_{\bar w}-\frac{A}{z-\bar z})\pi=0,&\\
\pi^\perp(\p_{\bar w} + \frac{|z|^2}{z-\bar z}\p_r - \frac{\bar
z}{z-\bar z} B)\pi=0. \eca
\end{equation}
Let
$$\bca
\cl= \p_{\bar w} -\frac{A}{z-\bar z},&\cr \cm=\p_{\bar w} +
\frac{|z|^2}{z-\bar z} \p_r -\frac{\bar z B}{z-\bar z}. \eca$$
Then
$$\cm-\cl= \frac{|z|^2}{z-\bar z}\p_r +\frac{1}{z-\bar z} (A-\bar z B).$$
So system \eqref{cu} is equivalent to
\begin{equation}\label{cv}
\bca
\pi^\perp(\p_{\bar w}- \frac{A}{z-\bar z})\pi=0,&\\
\pi^\perp(\p_r +\frac{1}{|z|^2}(A-\bar z B))\pi=0. \eca
\end{equation}
Thus we have shown

\bprop  \label{cz} Let $\psi, A, B$ be as in Proposition \ref{cy}.
Then $g_{z,\pi}\psi$ is an extended solution of the Ward equation
if and only if $\pi$ satisfies \eqref{cv}. \eprop

The first equation of \eqref{cv} has an interpretation in terms of
holomorphic subbundle.  To explain this, we first review some
notation of holomorphic vector bundles over $S^2$ (cf.
\cite{U1989}, \cite{Woo1989}).  A map $\eta$ defined and is smooth
on $S^2$ except on a finite subset $D$ is said to be of {\it pole
type\/} if at each $p_0\in D$ there exists a local complex
coordinate $(\co, w)$ of $S^2$ at $p_0$ with $w(p_0)=0$ such that
the map $\eta(w)= w^{-m} \eta_0(w)$ for all $w\in \co$, where $m$
is some positive integer and $\eta_0$ is smooth in a neighborhood
of $0$. Point $p_0$ is called a {\it pole\/} of $\eta$.

Given a smooth map $f:S^2\to U(n)$ and a constant $c$, let
$A=cf_{\bar w}f^{-1}$, then $\p_{\bar w} - A$ gives a holomorphic
structure on the trivial bundle $\underline{\C^n}= S^2\times \C^n$
over $S^2$.  A local section $\xi$ of $\underline{\C^n}$ is {\it
holomorphic in the complex structure $\p_{\bar w}-A$} if

\begin{equation}\label{cx}
\p_{\bar w}\xi -A\xi=0.
\end{equation}

  A {\it meromorphic section\/} of a holomorphic vector bundle  is a
section of pole type and is holomorphic away from the poles.  It
is known that the space of meromorphic sections of a rank $k$
holomorphic vector bundle $E$ over $S^2$ is of dimension $k$ over
the field ${\mathcal{R}}\/(S^2)$ of meromorphic functions on
$S^2$. In other words, there exist $k$ meromorphic sections
$\eta_1, \ldots, \eta_k$  such that if $\eta$ is a meromorphic
section of $E$ then there exist $f_1, \ldots,f_k$ in
${\mathcal{R}}(S^2)$ so that $\eta= \sum_{j=1}^k f_j \eta_j$. We
call such $\{\eta_1, \ldots, \eta_k\}$ a {\it meromorphic frame\/}
of $E$.

 The following  is known (cf. \cite{U1989,
Woo1989}):

\bprop \label{cw} Given a smooth map $\pi:S^2\to \Gr(k,\C^n)$, let
$\Pi$ denote the subbundle of $\underline{\C^n}$ whose fiber over
$p\in S^2$ is $\Im(\pi(p))$. Then the following two statements are
equivalent: \ben \item $\pi^\perp(\p_{\bar w} -A)\pi=0$, \item
$\Pi$ is a rank $k$ holomorphic subbundle of $\underline{\C^n}$
with respect to $\p_{\bar w}-A$. \een Moreover, if (1) or (2)
holds, then  there exist maps $\xi_1, \ldots, \xi_k:S^2\to \C^n$
of pole type so that \ben \item[(a)] $\xi_1(p), \ldots, \xi_k(p)$
span $\Im(\pi(p))$ for all $p\in S^2$ except at finitely many
points, \item[(b)] each $\xi_j$ is a solution of $\p_{\bar
w}\eta-A\eta=0$,
 \item[(c)]  every meromorphic section $\eta$ of $\Pi$ is of the form
$\sum_{j=1}^k f_j \xi_j$ for some $f_1, \ldots, f_k\in
{\mathcal{R}}(S^2)$. \een \eprop

As a consequence of the discussion above, we have

\bcor Let $\psi$ be an extended Ward soliton with
$A=(\l\psi_x-\psi_u)\psi^{-1}$ and $B=(\l\psi_v-\psi_x)\psi^{-1}$,
and $(p,q,r)$ the coordinate system on $\R^{2,1}$ defined by
\eqref{el}. Let $\pi: \R^{2,1}\to \Gr(k,\C^n)$ be a smooth map
that extends to $S^2\times \R$. Then the following statements are
equivalent:
 \ben
 \item
 $\pi$ is a solution of
\eqref{singBT}. \item For each fixed $r$, the subbundle $\Pi(r)$
associated to $\pi(\cdot,\cdot,r)$ is a holomorphic subbundle of
$\underline{\C^n}$ in the complex structure $\p_{\bar w}-
\frac{A}{z-\bar z}$ and satisfies
$\pi^\perp(\p_r+\frac{1}{|z|^2}(A-\bar z B))\pi =0$.
 \item
 There exist maps $\xi_1, \ldots, \xi_k: S^2 \times\R \to \C^n$
satisfying the following conditions:
 \ben
 \item  $\{\xi_1(\cdot,r), \ldots, \xi_k (\cdot,r)\}$ is a
meromorphic frame of $\Pi(r)$,
 \item
$\p_r\xi_j+\frac{1}{|z|^2}(A-\bar z B)\xi_j$ is a section of $\Pi$
of pole type  for all $1\leq j\leq k$.
 \een
 \een \ecor

\bs

 \section{Construction of all Ward solitons}

The goal of this section is to show that all Ward solitons can be
constructed by the methods given in sections 3, 5 and 6 (using
algebraic BTs, limiting method, and generalized algebraic BTs). By
Proposition \ref{bl}, it suffices to show that for any given $z\in
\C_+$ and $k\in \N$, we can construct all extended Ward solitons
of pole data $(z,k)$.

First we prove

\bthm \label{cf} If $\phi_k$ is an extended solution of the Ward
equation and $\phi_k=g_{z,\pi_k}\cdots g_{z,\pi_1}$ is the minimal
factorization, then the tails of $\phi_k$, $\phi_l=
g_{z,\pi_l}\cdots g_{z,\pi_1}$, $l=1,\cdots,k-1$, are also
extended solutions. \ethm

\begin{proof} We prove the theorem by induction on $l$. Since $\phi_k$ is an
extended solution,
\begin{equation}\label{cs}
\bca
A_k=(\l\p_x\phi_k -\p_u \phi_k)\phi_k^{-1}, &\\
B_k=(\l\p_v\phi_k -\p_x \phi_k)\phi_k^{-1}, \eca
\end{equation} are independent
of $\l$. We want to prove $g_{z,\pi_1}$ is an extended solution.
Compute the Laurent series expansion of $\phi_k$ at $\l=z$ to get
$$\phi_k=\I +\frac{z-\bar z}{\l-z}P_1 +\cdots +\left(
\frac{z-\bar z}{\l-z} \right)^k P_k,$$ where \begin{equation}
\label{cg} P_j =\sum_{k\ge i_1  >\cdots  >i_j \ge 1}
\pi_{i_1}^\perp \cdots \pi_{i_j}^\perp.\end{equation} Compute the
Laurent series expansion of both sides of
\begin{equation}\label{ci} \l\p_x\phi_k -\p_u \phi_k =A_k \phi_k
\end{equation} at $\l=z$, and compare the coefficients of
$(\l-z)^{-k}$ to see
\begin{equation}\label{cj}
(z\p_x -\p_u) P_k =A_kP_k,
\end{equation}
 where $P_k=\pi_k^\perp
\cdots\pi_1^\perp$. Multiply $\pi_1$ from right to both sides of
\eqref{cj} to see \begin{equation*}
 \pi_k^\perp \cdots \pi_2^\perp
((z\p_x -\p_u)\pi_1^\perp)\pi_1 =0.
 \end{equation*}
 But $\pi_1^\perp\pi_1=0$ implies
 $$L_z(\pi_1^\perp)\pi_1 + \pi_1^\perp L_z(\pi_1)=0,$$
where $L_z=z\p_x-\p_u$. So we have
\begin{equation}\label{ch}
\pi_k^\perp\cdots \pi_1^\perp L_z(\pi_1)=0.
\end{equation}
By assumption, $\phi_k =g_{z,\pi_k}\cdots g_{z,\pi_1}$ is a
minimal factorization, i.e.,
 $\Im\pi_j\cap\Im\pi_{j-1}^\perp=0$.  So by Lemma \ref{cr},
 $${\rm Ker\/}(\pi_k^\perp\cdots \pi_1^\perp) =\Im\pi_1.$$
 Hence equation \eqref{ch} implies that
 $$\Im(L_z(\pi_1))\subset \Im\pi_1.$$
 So  $$\pi_1^\perp(z\p_x \pi_1-\p_u \pi_1)=0.$$
 Use the second equation of \eqref{cs} and similar
argument to prove that $$\pi_1^\perp (z\p_v\pi_1 -\p_x\pi_1)=0.$$
The above two equalities imply that $\phi_1=g_{z,\pi_1}$ is an
extended solution of Ward equation.

Assume that $\phi_l$ is an extended solution, and let
$$A_l=(\l\p_x\phi_l -\p_u \phi_l)\phi_l^{-1}, \quad B_l=(\l\p_v\phi_l -\p_x
\phi_k)\phi_l^{-1}.$$
 We want to show that
$\phi_{l+1}$ is an extended solution too. Write
$\phi_k=\psi_{k-l}\phi_l$, where $\psi_{k-l}=g_{z,\pi_k}\cdots
g_{z,\pi_{l+1}}$.  Substitute $\phi_k= \psi_{k-l}\phi_l$ to
\eqref{ci} to get
 $$(\l\p_x
\psi_{k-l}-\p_u\psi_{k-l})\phi_l +\psi_{k-l} (\l\p_x
\phi_l-\p_u\phi_l) =A_k \psi_{k-l}\phi_l.$$ Multiply $\phi_l^{-1}$
from the right to both sides and use $\l\p_x \phi_l-\p_u\phi_l
=A_l \phi_l$ to get $$(\l\p_x \psi_{k-l}-\p_u\psi_{k-l})
+\psi_{k-l}A_l =A_k \psi_{k-l}.$$ Compute the Laurent series
expansion of both sides at $\l=z$, and compare the coefficients of
$(\l-z)^{-(k-l)}$ to get $$(z\p_x -\p_u)(\pi_k^\perp
\cdots\pi_{l+1}^\perp) +\pi_k^\perp \cdots \pi_{l+1}^\perp A_l
=A_k \pi_k^\perp \cdots \pi_{l+1}^\perp.$$ Multiply $\pi_{l+1}$
from right to both sides to get $$\pi_k^\perp \cdots
\pi_{l+1}^\perp (z\p_x \pi_{l+1} -\p_u \pi_{l+1} -A_l
\pi_{l+1})=0.$$ By Lemma \ref{cr}, we have
 $$\Im(z\p_x\pi_{l+1}-\p_u\pi_{l+1}-A_l \pi_{l+1})\subset {\rm
Ker\/}(\pi_k^\perp\cdots \pi_{l+1}^\perp)
  =\Im\pi_{l+1}.$$
Thus we have
$$\pi_{l+1}^\perp (z\p_x \pi_{l+1} -\p_u \pi_{l+1}
-A_l \pi_{l+1})=0.$$ Likewise, from $B_k$ we can obtain
$$\pi_{l+1}^\perp (z\p_v \pi_{l+1} -\p_x \pi_{l+1} -B_l
\pi_{l+1})=0.$$ The above two equalities imply that $\pi_{l+1}$ is
a solution of the analytic BT \eqref{singBT} with
$(A,B)=(A_l,B_l)$, hence $\phi_{l+1} =g_{z,\pi_{l+1}}\phi_l$ is an
extended solution. Thus we complete the proof by induction.
\end{proof}

Theorem~\ref{cf} tells us that any extended solution of pole data
$(z,k)$ is obtained by solving the analytic BT \eqref{singBT} of
an extended $1$-soliton, then of an extended $2$-soliton, ... etc.

 The following two Lemmas prove that to solve the analytic BT
\eqref{singBT}, which is a system of non-linear equations,  it
suffices to solve certain first order linear system.

\begin{lem} \label{cn}
Let $\psi$ be an extended solution of the Ward equation, and
$A=(\l\psi_x-\psi_u)\psi^{-1}$, and
$B=(\l\psi_v-\psi_x)\psi^{-1}$. If $\pi$ is a local solution of
\eqref{singBT}, then there exists a local $\cm^0_{n\times
r}(\C)$-valued smooth map $V$ so that columns of $V$ spans
$\Im\pi$ and $V$ satisfies

 \begin{equation}\label{cm}
\bca
    zV_x-V_u-AV=0, \\
    zV_v-V_x-BV=0.
\eca
\end{equation}
 Conversely, if $V$ is a solution of \eqref{cm}, then the projection
$\pi= V(V^*V)^{-1}V^*$ is a solution of \eqref{singBT}.
  \end{lem}

  \begin{proof}
  Choose  a local $\cm^0_{n\times r}(\C)$-valued smooth map $V$ such
that columns of $V$ span $\Im\pi$.  So $\pi=V(V^*V)^{-1}V^*$.
Substitute this into \eqref{singBT} to see
\begin{equation} \label{cl} \left\{%
\begin{array}{l}
    \pi^\perp (zV_x-V_u-AV)(V^*V)^{-1}V^*=0, \\
    \pi^\perp (zV_v-V_x-BV)(V^*V)^{-1}V^*=0. \\
\end{array}%
\right. \end{equation}
  Multiply $V$ from the right to see
$$  \left\{%
\begin{array}{l}
    \pi^\perp (zV_x-V_u-AV)=0, \\
    \pi^\perp (zV_v-V_x-BV)=0. \\
\end{array}%
\right.$$
  This implies that $\Im(zV_x-V_u-AV)$ and
$\Im(zV_v-V_x-BV)$ lie in $\Im\pi$. So
\begin{equation} \label{bw}
  \bca
    zV_x-V_u-AV=Vh, \\
    zV_v-V_x-BV=Vk,
  \eca
  \end{equation}
   for some $r\times r$ matrix-valued maps $h,k$.

  Claim that there exists a smooth $GL(r,\C)$-valued map $\phi$ so that
  $\ti V= V\phi$ satisfies \eqref{cm}.  To see this, let
  $$L_1= z\p_x-\p_u, \quad L_2= z\p_v-\p_x.$$
  Since $L_1, L_2$ are constant coefficient linear operators, they
commute.  A direct computation shows that $V\phi$ satisfies
\eqref{cm} if and only if
  $\phi$ satisfies
  \begin{equation}\label{bv}
  \bca
  L_1\phi=- h\phi,\\
  L_2 \phi= -k\phi.
  \eca
  \end{equation}
Equation \eqref{bv} is solvable if and only if $h, k$ satisfy
 \begin{equation}\label{by}
- L_2(h)+hk=- L_1(k)+ kh.
 \end{equation}
  This condition for $h,k$ comes from equating $L_1L_2\phi= L_2L_1\phi$.

 Write \eqref{bw} in terms of $L_1, L_2$ to get
 $$\bca L_1(V)= AV+Vh,\\ L_2(V)= BV+ Vk.\eca$$
 Since $L_1, L_2$ commute, we have
 \begin{align*}
 L_2L_1(V)&= L_2(A) V+ A L_2(V) + L_2(V)h + VL_2(h) \\
 &= (L_2(A) + AB)V + AVk + (BV+Vk)h + VL_2(h)\\
 = L_1L_2(V) &= L_1(B)V + B L_1(V) + L_1(V) k + V L_1(k)\\
 &= (L_1(B) + BA)V + BVh+ (AV+ Vh)k+ VL_1(k).
 \end{align*}
 Hence
 \begin{equation}\label{bx}
 (L_2(A)-L_1(B)+[A,B])V= V(L_1(k)-L_2(h) +[h,k]).
 \end{equation}
 Since $\psi$ is an extended Ward map, $A, B$ satisfies
\eqref{eq:auxiliary}, which implies that
 \begin{equation}\label{bz}
 L_2(A)-L_1(B)+[A,B]=0.
 \end{equation}
By \eqref{bx} and \eqref{bz}, $h, k$ satisfy \eqref{by}. Thus we
can find local smooth solution $\phi$ for \eqref{bv} and $\ti V=
V\phi$ solves \eqref{cm}.

The converse is clearly true.
 \end{proof}

   It is easy to see that solutions of \eqref{cm} are unique up to
$V\to VH$, where $H$ is a map from $\R^{2,1}$ to $GL(r,\C)$
satisfying
$$\left\{%
\begin{array}{l}
    zH_x-H_u=0, \\
    zH_v-H_x=0, \\
\end{array}%
\right.$$ i.e. $H$ depends on $w=x+zu+z^{-1}v$ only.

If we can construct local fundamental solutions to the linear
system \eqref{cm}, then we can obtain all local extended solutions
of the form $g_{z,\pi}\psi$.  But in order to construct global
Ward maps satisfying the boundary condition \eqref{bc1}, we need
to construct fundamental solutions that are of pole type on each
$w$-plane.

If $\psi$ is an extended Ward soliton, holomorphic and
non-degenerate at $\l=z$, then $\psi(z)$ itself is a fundamental
solution to \eqref{cm}. So we have:

\bprop \label{df}
  Let $\psi$ be an extended Ward soliton, and $A=
(\l\psi_x-\psi_u)\psi^{-1}$ and $B=(\l \psi_v-\psi_x)\psi^{-1}$.
If $\psi(x,u,v,\l)$ is holomorphic and non-degenerate at $\l=z$,
then $\psi(\cdots, z)$ is a fundamental solution of \eqref{cm}.
Consequently, the columns of $\psi(\cdots, z)$ form a holomorphic
frame of the trivial $\C^n$-bundle over $S^2$ with respect to the
complex structure $\p_{\bar w}-\frac{A}{z-\bar z}$ on each
$w$-plane.
  \eprop

\begin{proof}
Since the extended solution $\psi$ is holomorphic and
non-degenerate at $\l=z$, \eqref{eq:auxiliary} implies that
$$\bca
z\psi_x(z) -\psi_u(z)-A\psi(z)=0,&\\
z\psi_v(z) -\psi_u(z)-B\psi(z)=0. \eca$$ Hence the columns of
$\psi(\cdots, z)$ are smooth solutions of \eqref{cm}. So on each
$w$-plane, columns of $\psi(z)$ form a holomorphic frame of the
trivial $\C^n$-bundle over $S^2$ with respect to the complex
structure $\p_{\bar w} -\frac{A}{z-\bar z}$.
\end{proof}

The above Proposition implies that Theorem~\ref{thm:BT} gives all
extended solutions of the form $g_{z,\pi}\psi$ when $\psi$ is
holomorphic and non-degenerate at $\l=z$.

\ms

If $\psi$ is an extended Ward soliton and has a pole at $\l=z$,
then we will show below that the limiting method used in section 5
gives fundamental solutions of \eqref{cm} that are of pole type on
each $w$-plane.

We need two lemmas first, and their proofs are straight forward.

 \blem\label{de}
Let $b_\e$, $\eta_\e$ and $\psi_\e$ be  maps from
$\R^{2,1}\times\Omega$ to $gl(n,\C)$, where $\Omega$ is an open
subset of $0$ in $\C$. Suppose $b_\e= \sum_{j=0}^\infty c_j\e^j$,
$\psi_\e= \sum_{j=0}^k P_j\e^{-j}$, and $\eta_\e= \psi_\e(b_\e)$.
If $\psi_\e(b_\e)$ is smooth at $\e=0$, i.e., $\sum_{j-i=m} P_j
c_i =0$ for all $m>0$, then
  \ben
\item $\eta:= \lim_{\e\to 0} \eta_\e = \sum_{j=0}^k P_j c_j$,
\item $\bca
z\eta_x- \eta_u= \lim_{\e\to 0} (z+\e)(\eta_\e)_x- (\eta_\e)_u, &\\
z\eta_v-\eta_x= \lim_{\e\to 0} (z+\e)(\eta_\e)_v- (\eta_\e)_x.
\eca$ \een \elem

\blem\label{cp} Let $w_\e= x+(z+\e)u+(z+\e)^{-1}v$, and $b_{i\e}=
\sum_{j=0}^{k_i} \e^j a_{ij}(w_\e)$ for $i=1,2$.  Let $\phi$ be an
extended Ward map, and $f_1, f_2$ meromorphic functions from $\C$
to $\C$. If $\lim_{\e\to 0}\phi(z+\e) (b_{i\e})=\eta_i$ exists for
$i=1, 2$, then
$$\lim_{\e\to 0} \phi(z+\e)(f_1(w_\e)b_{1\e}+f_2(w_\e) b_{2\e})=
f_1\eta_1+ f_2 \eta_2.$$ \elem

\begin{proof}
The Lemma follows from
\begin{align*}
&\phi(z+\e)(f_1(w_\e)b_{1\e}+ f_2(w_\e)b_{2\e})=
f_1(w_\e)\phi(z+\e)(b_{1\e}) + f_2(w_\e)\phi(z+\e)(b_{2\e})\\
&\quad \to f_1(w)\eta_1+ f_2(w)\eta_2 \quad {\rm as\ \ }\e\to 0.
\end{align*}
\end{proof}

\bthm\label{di} Let $\phi_k$ be an extended Ward soliton with pole
data $(z,k)$, and
$$A_k= (\l(\phi_k)_x -(\phi_k)_u)\phi_k^{-1}, \quad
B_k=(\l(\phi_k)_v- (\phi_k)_x)\phi_k^{-1}.$$ Then
  \ben
  \item
$\phi_k$ can be constructed using algebraic BT and the limiting
method given in section 5,
  \item we can use algebraic BT and the limiting method to
construct a fundamental solution $\eta$ of  \eqref{cm} with
$A=A_k, B=B_k$ and the entries of $\eta$ are  rational functions
in $x,y$ and $t$.
  \een \ethm

\bproof We prove the Theorem by induction on $k$.  For $k=1$, (1)
is obvious. For (2), we first choose $\C^n$-valued rational maps
$b_{n_1+1}, \ldots, b_n$ such that  $\pi_1^\perp(b_{n_1+1}),
\ldots , \pi_1^\perp(b_n)$ span $\Im \pi_1^\perp$, where
$n_1=\rk(\pi_1)$. Let
$$u_{j, \e}= g_{z,\pi_1}(z+\e)(\e b_j(w_\e)).$$
   When $|\e|>0$ is
small, $g_{z,\pi_1}$ is holomorphic and non-degenerate at
$\l=z+\e$. Hence by Proposition \ref{df}, we have
$$\bca
(z+\e)(u_{j,\e})_x - (u_{j,\e})_u -A_1 u_{j,\e}=0,&\\
(z+\e)(u_{j,\e})_v-(u_{j,\e})_v- B_1 u_{j,\e}=0. \eca$$
     A direct computation implies that
$$u_j=\lim_{\e\to 0} u_{j,\e}= (z-\bar z)\pi^\perp_1(b_j(w)).$$
      As $\e\to
0$, Lemma \ref{de} implies that $u_j$ satisfies
$$\bca
z(u_j)_x-(u_j)_u-A_1 u_j=0,&\\
z(u_j)_v-(u_j)_x-B_1 u_j=0, \eca$$ for $n_1+1\leq j\leq n$.  So
$u_{n_1+1}, \ldots, u_n$ are rational maps in $x,y,t$ and are
linearly independent solutions of \eqref{cm} with $A=A_1$ and
$B=B_1$.

We claim that by choosing a sequence of rational functions
carefully, we can construct the rest linearly independent
solutions of \eqref{cm}. Let $a_1, \ldots, a_{n_1}$ be
$\C^n$-valued rational maps that span $\Im\pi_1$ except at
finitely many points. Use formulas in section 5 to get
$$ g_{z,\pi_1}(z+\e)(a_j(w_\e))= (\I+\frac{z-\bar
z}{\e}\pi_1^\perp)(a_j(w)+(u-z^{-2}v)a_j'(w)\e +O(\e^2)).$$
      So
$$\eta_j:=\lim_{\e\to 0} g_{z,\pi_1}(z+\e)(a_j(w_\e))=  a_j(w)
+(z-\bar z)(u-z^{-2}v) \pi_1^\perp(a_j'(w))$$
 is rational in $x,y$ and $t$.
      By Proposition \ref{df}
and Lemma \ref{de}, $\eta_j$ is a solution of \eqref{cm} with
$A=A_1$ and $B=B_1$. Because $u_\a$ with $n_1+1\leq \a\leq n$ span
$\Im\pi_1^\perp$ and $\pi_1(\eta_1), \ldots, \pi_1(\eta_{n_1})$
span $\Im\pi_1$, $\{ \eta_1, \ldots, \eta_{n_1}, u_{n_1+1}, \ldots
u_n \}$ form a fundamental solution of \eqref{cm}. We have proved
that the entries of $\eta_i$ and $u_j$ are rational in $x,y$ and
$t$. This proves the claim and the Theorem for $k=1$.

Suppose the Theorem is true for $k$.   We want to prove that (1)
holds for $k+1$.   We may assume that
$$\phi_{k+1}= g_{z,\pi_{k+1}}\cdots g_{z,\pi_1}$$
 is the minimal factorization.  Let $n_j=\rk(\pi_j)$.  By Theorem
\ref{cf}, $\phi_k= g_{z,\pi_{k}}\cdots g_{z,\pi_1}$ is also an
extended solution and is a minimal factorization. By induction
hypothesis, we can construct a fundamental solution $\eta=(\eta_1,
\ldots, \eta_n)$ rational in $x, y, t$ for \eqref{cm} with $A=A_k$
and $B=B_k$.  Since $\phi_k$ and $g_{z,\pi_{k+1}}\phi_k$ are
extended solutions, by Proposition \ref{cy}, $\pi_{k+1}$ is a
solution of \eqref{singBT}.  By Proposition \ref{cw}, there exist
maps $\xi_1, \ldots, \xi_{n_{k+1}}$ of pole type that span
$\Im\pi_{k+1}$ and satisfy \eqref{cm} with $A=A_k$ and $B=B_k$.
But $\eta$ is a fundamental solution of \eqref{cm} over the field
${\mathcal R}(S^2)$.  So there exists a rational map $h=(h_{ij})$
from $\C$ to  $\cm^0_{n\times n_{k+1}}(\C)$ such that $\xi_j=
\sum_{i=1}^n h_{ij} \eta_i$ for $1\le j\le n_{k+1}$. By induction
hypothesis, each $\eta_i$ is constructed by the limiting method.
It follows from Lemma \ref{cp} that $\xi_j$ can be constructed by
the limiting method. This proves (1) for $k+1$.

 To prove (2), let $g=g_{z,\pi_{k+1}}$.  Note that $\phi_{k+1}=
g\phi_k$ implies
 $$L_\l(g)g^{-1}+ gL_\l(\phi_k)\phi_k^{-1}g^{-1}= A_{k+1},$$
 where $L_\l(\xi)= \l\xi_x-\xi_u$.  But
$L_\l(\phi_k)\phi_k^{-1}=A_k$. So we have
 $$L_\l(g) + gA_k = A_{k+1} g.$$
Equate the residue of the above equation at $\l=z$ to get
\begin{equation}\label{dc}
L_z(\pi_{k+1}^\perp) +\pi_{k+1}^\perp A_k= A_{k+1}\pi_{k+1}^\perp.
\end{equation}
We have $L_z(\eta)= A_k\eta$. Set $W= \pi_{k+1}^\perp \eta$, where
$\eta$ is a fundamental solution of \eqref{cm} for $\phi_k$.  We
want to show that $W$ satisfies \eqref{cm} with $A=A_{k+1}$ and
$B=B_{k+1}$.  To see this, we compute
\begin{align*}
L_zW- A_{k+1}W &= L_z(\pi^\perp_{k+1}\eta)
-A_{k+1}\pi_{k+1}^\perp \eta \\
 &= L_z(\pi_{k+1}^\perp)\eta+ \pi_{k+1}^\perp L_z(\eta) -
A_{k+1}\pi_{k+1}^\perp \eta\\
 &= (L_z(\pi_{k+1}^\perp)+ \pi_{k+1}^\perp A_k- A_{k+1}\pi_{k+1}^\perp)\eta,
 \end{align*}
 which is zero by \eqref{dc}. Similar argument implies that $zW_v-W_x
-B_{k+1}W=0$. This proves the claim.  We may assume that the
columns $\zeta_1, \ldots, \zeta_{n-n_{k+1}}$ of $W$ are linearly
independent. So these columns are linearly independent solutions
of \eqref{cm} with $A=A_{k+1}$ and $B=B_{k+1}$. Since $\pi_{k+1}$
is constructed by limiting method, so are the $\zeta_j$'s.  It
remains to construct $n_{k+1}$ other linearly independent
solutions of \eqref{cm}. We have proved (1) for $k+1$. So there
exist
 $$b_{j,\e}= \sum_{i=0}^k \e^i a_{ji}(w_\e)$$ with rational maps
$a_{ji}$ from $\C$ to $\C^n$ so that
 $$\phi_k(z+\e)(b_{j,\e})= \xi_j + \e y_j + O(\e^2),$$
 and $\xi_1, \ldots, \xi_{n_{k+1}}$ span $\Im\pi_{k+1}$.
 The limit of
 $$\phi_{k+1}(z+\e)(b_{j,\e})= (\I+\frac{z-\bar z}{\e}
\pi_{k+1}^\perp)(\xi_j +\e y_j+ O(\e^2))$$
 as $\e\to 0$ is
 $$\eta_j= \xi_j +(z-\bar z)\pi_{k+1}^\perp(y_j),
 \quad 1\leq j\leq n_{k+1}.$$
 By Lemma \ref{de}, these $\eta_j$'s are solutions of
\eqref{cm}.  So $$\eta_1, \ldots, \eta_{n_{k+1}}, \zeta_1, \ldots,
\zeta_{n-n_{k+1}}$$ form a fundamental solution of \eqref{cm} with
$A=A_{k+1}$ and $B=B_{k+1}$.  The above arguments also prove that
all entries of  $\eta_i$ and $\zeta_j$ are rational in $x,y$ and
$t$. This completes the proof of the Theorem.
 \eproof

 As a consequence of Proposition \ref{bl} and Theorem \ref{di}, we have

 \bcor
 All Ward solitons can be constructed by algebraic B\"acklund
transformations, the limiting method and generalized B\"acklund
transformations in sections 3, 5 and 6. Moreover, the entries of
Ward solitons are rational functions in $x, y$ and $t$.
 \ecor

\bs

\section{An explicit construction of unitons}

We have proved in the last section that all Ward solitons of pole
type $(i, k)$ can be constructed by the limiting method of section
5.  Since $U(n)$-unitons are stationary Ward solitons with pole
type $(i, k)$ for some $k<n$, we can use our method to construct
unitons. In this section, we write down the conditions on the
sequence of rational maps so that the limiting Ward solitons are
independent of $t$.  We then give explicit formulas for unitons
arising from finite sequence of $\C^n$-valued rational maps.  We
note that Wood \cite{Woo1989} and Burstall-Guest \cite{BurGue97}
also gave algorithms to construct unitons. Our
construction is somewhat different from theirs.

It is proved by Uhlenbeck that every $k$-uniton has a unique
extended solution of the form $\phi_k=g_{i,\pi_k}\cdots
g_{i,\pi_1}$ with the property that the span of
$$\{\pi_i\pi_{i-1}\cdots \pi_1(w)\n w\in S^2\}$$
is $\C^n$ for each $1\leq i\leq k$.  Moreover, $\rk\pi_1>\rk\pi_2
>\cdots  >\rk\pi_k$.  We will prove later that minimal factorizations
and explicit constructions also give the same condition on the
ranks of projections. But the condition we have on the $\pi_i's$
is $\Im\pi_{i+1}\cap \Im\pi_i^\perp=0$ for all $1\leq i\leq k-1$.

We use the same notation as in section 5, and assume that the pole
$z=i$. Let $a_0, a_1, \ldots$ be a sequence of $\C^n$-valued
meromorphic functions on $S^2$,
$$w_\e= x+(i+\e)u + (i+\e)^{-1}v, \quad w= x+iy,$$
\begin{align*}
b_{j,\e}&= a_0(w_\e)+ \e a_1(w_\e) + \cdots + \e^j a_j(w_\e),\\
&= c_0+ c_1 \e + \cdots + c_j\e^j +O(\e^{j+1}),
\end{align*}
where $c_j$'s are defined by \eqref{ao} with $z=i$.  So
\begin{eqnarray*} c_0 &=& a_0(w), \\
c_1 &=& ta^\prime_0(w) +a_1(w), \\
c_2 &=& \frac{t^2}{2}a^{\prime\prime}_0 (w) +iva^\prime_0(w)
+ta^\prime_1(w) +a_2(w),\\
c_3 &=& \frac{t^3}{3!} a^{\prime\prime\prime}_0 (w)
+ivta^{\prime\prime}_0 (w) -va^\prime_0 (w)
+\frac{t^2}{2}a^{\prime\prime}_1 (w) +iva^\prime_1(w)\\
& & +ta^\prime_2(w) +a_3(w),\\
c_4 &=& \frac{t^4}{4!} a^{(4)}_0 (w) +\frac{i}{2}vt^2
a^{\prime\prime\prime}_0 (w) -\frac{a^{\prime\prime}_0(w)}{2} (v^2
+2vt) -iv a^\prime_0 (w) \\ & & +\frac{t^3}{3!}
a^{\prime\prime\prime}_1 (w) +ivta^{\prime\prime}_1 (w)
-va^\prime_1 (w) +\frac{t^2}{2}a^{\prime\prime}_2 (w)
+iva^\prime_2(w)\\ && +ta^\prime_3(w) +a_4(w),\\ & & \cdots
\end{eqnarray*}

We want to write down the conditions that the limiting Ward
soliton is independent of $t$.  Suppose
$$\phi_k= g_{i, \pi_k}g_{i,\pi_{k-1}}\cdots g_{i,\pi_1}$$ is a
minimal factorization and is an extended Ward $k$-soliton obtained
by the limiting method of section 5. Let $V_j=\Im\pi_j$, and
\begin{align*}
P_{ji} &=\sum_{j\geq k_i> \cdots  > k_1\geq 1}
\pi_{k_i}^\perp\cdots
\pi_{k_1}^\perp,\\
P_{j,0}&=\I.
\end{align*}
 Use the computation of section 5 to conclude that:

\ben \item If $\rk\pi_2=1$, then $\phi_2$ is independent of $t$ if
and only if
\begin{align*}
&a_0, a_0'\in V_1,\\
& v_2= a_0+ 2i\pi_1^\perp(a_1).
\end{align*}
\item If $\rk\pi_3=1$, then $\phi_3$ is independent of $t$ if and
only if
$$v_3= v_2+ 2i\pi_2^\perp(c_1+ 2i\pi_1^\perp(c_2))$$
is independent of $t$.  So the coefficients of $t, t^2$ must be
zero, which gives
\begin{equation*}
\bca
\pi_2^\perp\pi_1^\perp(a_0'')=0,&\\
\pi_2^\perp(a_0'+ 2i\pi_1^\perp(\frac{i}{2} a_0'+ a_1'))=0. \eca
\end{equation*}
Since $\phi_2$ is the minimal factorization, by Proposition
\ref{cr} (2), the first equation implies that $a_0''\in V_1$.  We
have $a_0'\in V_1$.  So the condition for $\phi_3$ to be
independent of $t$ is
\begin{equation}\label{ei2}
\bca
a_0, a_0', a_0''\in V_1,&\\
a_0+2i\pi_1^\perp(a_1), a_0'+ 2i\pi_1^\perp(a_1')\in V_2, \eca
\end{equation}
and
$$v_3= a_0+ 2i(\pi_1^\perp+\pi_2^\perp)(a_1)+
(2i)^2\pi_2^\perp\pi_1^\perp( a_2).$$

\item If $\rk\pi_4=1$, then $\phi_4$ is independent of $t$ if and
only if all coefficients of $t, t^2, t^3$ in
$$v_4= v_3+ 2i\pi_3^\perp(c_1)+(2i)^2(\pi_3^\perp \pi_1^\perp+
\pi_3^\perp\pi_2^\perp) (c_2) +
(2i)^3\pi_3^\perp\pi_2^\perp\pi_1^\perp(c_3)$$ are zero.  So we
get
\begin{equation}\label{ei3}
\bca
a_0, a_0', a_0'', a_0'''\in V_1,&\\
a_0^{(i)}+ 2i\pi_1^\perp(a_1^{(i)})\in V_2,& 0\leq i\leq 2\\
a_0^{(i)}+ 2i(\pi_1^\perp+\pi_2^\perp)(a_1^{(i)}) +
(2i)^2\pi_2^\perp\pi_1^\perp(a_2^{(i)})\in V_3, & i=0,1 \eca
\end{equation}
and
\begin{equation}\label{br}
v_4= a_0+ (2i)P_{31}(a_1)+(2i)^2P_{32}(a_2)+ (2i)^3P_{33}(a_3).
\end{equation}
\item By induction, if $\rk\pi_k=1$, then $\phi_k$ is independent
of $t$ if and only if
\begin{equation*}
\bca
a_0, \ldots, a_0^{(k-1)}\in V_1,\\
v_2, D^{(1)} v_2,\ldots, D^{(k-2)} v_2\in V_2,\\
\cdots,\\
v_{k-1}, D^{(1)} v_{k-1}\in V_{k-1}, \eca
\end{equation*}
where
$$D^{(i)} v_\ell= \sum_{j=0}^{\ell-1} (2i)^j P_{\ell-1,j}(a_j^{(i)}),$$
and $a^{(i)}= \frac{d^ia}{dw^i}$. Thus
$$v_k= \sum_{j=0}^{k-1} (2i)^j P_{k-1,j}(a_j).$$

\een

The computation for the case when $\rk\pi_k\geq 2$ is similar.  In
fact, we get

\bthm \label{ds} Suppose $\phi_k$ is an extended uniton and
$\phi_k=g_{i,\pi_k}\cdots g_{i,\pi_1}$ is a minimal factorization.
Let $m_i =\rk \pi_i$, and
$$P_{m,j}=\sum_{m\geq i_1> \cdots  > i_j\geq 1} \pi_{i_1}^\perp\cdots
\pi_{i_j}^\perp$$ for $1\leq m\leq k$. Then there exists a
partition $(r_1, \ldots, r_s)$ of $m_k$,  (i.e., $r_j>0$ and
$\sum_{j=1}^s r_j= m_k$) and $\C^n$-valued  rational maps
$a_{i,0}, \ldots, a_{i,k-1}$ for $1\leq i\leq s$ such that \ben
\item $a_{i,0}^{(j_i)}\in \Im\pi_1$ for all $0\leq j_i\leq
r_i+k-2$, \item $D^{(j_i)}v_{ip}\in \Im\pi_p$ for all $0\leq
j_i\leq r_i+k-p-1$ and $1\leq p\leq k-1$, \item $\{v_{ik},
D^{(1)}v_{ik}, \ldots, D^{(r_i-1)}v_{ik}\n 1\leq i\leq s\}$ spans
$\Im\pi_k$ and $D^{(r_i)}v_{ik}\not\in \Im\pi_k$, \een where
$$v_{im}=\sum_{j=0}^{m-1} (2i)^j P_{m-1, j}(a_{i,j}), \quad
D^{(\ell)} v_{im}= \sum_{j=0}^{m-1} (2i)^jP_{m-1,
j}a_{i,j}^{(\ell)}.$$ \ethm

Let $\phi_k= g_{i,\pi_k}\cdots g_{i,\pi_1}$ be a minimal
factorization and an extended uniton.  $(m_1, \ldots, m_k)$ is
called the {\it rank data\/} of $\phi_k$, where  $m_j=\rk(\pi_j)$.
By Proposition \ref{cr}, $m_1\geq \cdots \geq m_k$.  We will prove
below that $m_i$ are strictly decreasing.  To do this, we first
note that we may assume $V_{11}=\Im\pi_1\cap \C^n=0$, where $\C^n$
means constant maps from $S^2$ to $\C^n$. If not, then
$g_{i,\pi_1}=\frac{\l -i}{\l +i} g_{i,\tau_1}g_{i,\tau_2}$, where
$\tau_1$ and $\tau_2$ are projections onto $V_{11}^\perp\cap
\Im\pi_1$ and $V_{11}$ respectively.  Since $\tau_2$ is a constant
projection, the harmonic maps corresponding to $\phi_k$ and to
$g_{i, \pi_k}\cdots g_{i,\pi_2} g_{i, \tau_1}$ are only differed
by the left multiplication of a constant element in $U(n)$.  So we
may assume that $V_{11}=0$.

\bprop \label{dm} Under the same assumption as in Theorem
\ref{ds}, if $\Im\pi_1\cap \C^n=0$,  then $m_1>\cdots  > m_k$.
\eprop

\begin{proof}  It follows from Proposition \ref{cr} that $m_1\geq
m_2\geq \cdots \geq m_k$. We prove the Proposition by induction on
$k$. For $k=2$, if the Proposition is not true, then $m_1=m_2$. By
Theorem \ref{ds}, there exists a partition $(r_1, \ldots, r_s)$ of
$m_2$ and $\C^n$-valued rational maps $a_{i,0}, a_{i,1}$ for
$1\leq i\leq s$ so that \ben \item[(i)] $a_{i,0}^{(j_i)}\in
\Im\pi_1$ for $0\leq j_i\leq r_i$, \item[(ii)]
$\{D^{(j_i)}v_{i,2}\n  1\leq i\leq s, \  0\leq j_i\leq r_i-1\}$
spans $\Im\pi_2$ a.e., where
 $$v_{i2}= a_{i,0}+2i\pi_1^\perp(a_{i,1}), \quad D^{(i)}v_{\ell,2}=
a_{\ell,0}^{(i)}+ 2i\pi_1^\perp(a_{\ell, 1}^{(j)}),$$
 \item[(iii)] $D^{(r_i)}v_{i,2}\not\in \Im\pi_2$.
 \een
  The definition of minimal factorization implies
$\Im\pi_1^\perp\cap\Im\pi_2=0$.
 So the rank of $\{a_{i,0}^{(j_i)}\n 1\leq i\leq s, 0\leq j_i\leq
r_i-1\}$ is $m_2$, which is equal to $m_1=\rk\pi_1$.  But
$D^{(r_i-1)}v_{i2}\in \Im\pi_2$ implies that $a_{i,0}^{(r_i)}\in
\Im\pi_1$. It follows from a direct computation  that $\p_w V=hV$
for some meromorphic function $h$, where
$$V= a_{1,0}\wedge\cdots \wedge a_{1,0}^{(r_1-1)}\wedge \cdots \wedge
a_{s,0}\wedge\cdots \wedge a_{s,0}^{(r_s-1)}.$$  Define
$$f(w)= \exp(-\int h(w)dw)$$ locally.
Then $\p_w (fV)=0$.
 But $fV$ is meromorphic.  Hence $fV$ is locally constant, which implies
that $\pi_1$ is a constant Hermitian projection.  Thus
$\Im\pi_1\subset \C^n$. In particular, $\Im\pi_1\cap \C^n\not=0$,
a contradiction. So $m_1>m_2$.

 Suppose the Proposition is true for $k$, and
 $$\phi_{k+1}= g_{i,\pi_{k+1}}\cdots g_{i,\pi_1}$$
 satisfies the assumption of the Proposition.  It follows from
Theorem \ref{cf} that
  $$\phi_k= g_{i,\pi_k}\cdots g_{i,\pi_1}$$
is also an extended uniton. It is easy to check that $\phi_k$ also
satisfies the conditions of the Proposition.  So by the induction
hypothesis, $m_1>\cdots  > m_k$. We already have $m_k\geq m_{k+1}$
by Proposition \ref{cr}. We will show that  $m_k=m_{k+1}$ gives a
contradiction next. By Theorem \ref{ds}, there exist a partition
$(r_1, \ldots, r_s)$ of $m_{k+1}$ and $a_{i, j}$'s so that
$\{D^{(j_i)}v_{i,k+1}\n 1\leq i\leq s, 0\leq j_i\leq r_i-1\}$ form
a basis of $\Im\pi_{k+1}$. Use
$$P_{k,j}= P_{k-1, j} + \pi_k^\perp P_{k-1, j-1}$$
to see that $D^{(\ell)}v_{i,k+1}=D^{(\ell)} v_{i,k}  +\xi_{i,k,
\ell}$ for some $\xi_{i,k,\ell}$ in $\Im\pi_k^\perp$.  Since
$\Im\pi_{k+1}\cap \Im\pi_k^\perp=0$,
$$\{D^{(\ell_i)} v_{i, k}\n 1\leq i\leq s, 0\leq \ell_i\leq r_i-1\}$$
has rank $m_{k+1}$, which is equal to $m_k$.  But
$D^{(r_i-1)}v_{i,k+1}\in \Im\pi_{k+1}$ implies that
$D^{(r_i)}v_{i,k}\in\Im\pi_k$.  Induction hypothesis says that
this can not happen, a contradiction.
\end{proof}

We give some examples to demonstrate how to  write down unitons
from rational maps.

\beg An extended $4$-uniton $\phi_4=g_{i,\pi_4}\cdots g_{i,\pi_1}$
in $U(5)$ with rank data  $(4,3,2,1)$ is given by $\C^5$-valued
rational maps $a_0, a_1, a_2, a_3$ such that \ben \item $\Im\pi_1$
is spanned by $a_0, \ldots, a_0^{(3)}$, \item $\Im\pi_2$ is
spanned by $a_0^{(j)}+2i\pi_1^\perp(a_1^{(j)})$ with $0\leq j\leq
2$, \item $\Im\pi_3$ is spanned by $a_0^{(j)}+
2iP_{21}(a_1^{(j)})+ (2i)^2P_{22}(a_2^{(j)})$ with $j=0, 1$, \item
$\Im\pi_4$ is spanned by $v_4$ defined by \eqref{br}. \een Note
that $a_0$ should be chosen so that $a_0,a_0',\cdots,a_0^{(4)}$
are linearly independent a.e., otherwise $\pi_1$ is constant and
contradicts $\Im \pi_1\cap \C^n =0$. \eeg

\beg A $3$-uniton $\phi_3=g_{i,\pi_3}g_{i,\pi_2}g_{i,\pi_1}$ in
$U(5)$ with rank data $(4, 2,1)$ is given by $\C^5$-valued
rational maps $a_0, b_0, a_1, a_2$ \ben \item $\Im\pi_1$ is
spanned by $a_0, a_0', a_0'', b_0$, \item $\Im\pi_2$ is spanned by
$a_0+ 2i\pi_1^\perp(a_1)$, $a_0'+ 2i\pi_1^\perp(a_1')$, \item
$\Im\pi_3$ is spanned by $\sum_{j=0}^2 (2i)^jP_{2j}(a_j)$. \een
\eeg

\beg $3$-unitons in $U(6)$ with rank data $(5, 3,1)$. \par

\ms \ni \underline{Case 1}. Choose $\C^6$-valued rational maps
$a_0, a_1, a_2, b_0, b_1$ on $\C$ such that \ben \item $a_0, a_0',
a_0'', b_0, b_0'$ are linearly independent a.e., and their span
intersects the space $\C^n$ of constant maps from $\C$ to $\C^n$
only at $0$, \item $v_2=a_0+2i\pi_1^\perp(a_1)$, $D^{(1)}v_2=
a_0'+ 2i\pi_1^\perp(a_1')$, $\hat v_2=b_0+ 2i\pi_1^\perp(b_1)$ are
linearly independent a.e., \item $v_3=a_0+2iP_{21}(a_1)+
(2i)^2P_{22}(a_2)$ is not zero a.e.. \een Let $\pi_1$, $\pi_2$,
and $\pi_3$ be the projections of $\C^6$ onto the span of $a_0,
a_0', a_0'', b_0, b_0'$, the span of $v_2, D^{(1)}v_2, \hat v_2$,
and $\C v_3$ respectively.  Then $\phi_3= g_{i, \pi_3}g_{i,\pi_2}
g_{i,\pi_1}$ is an extended solution of a $3$-uniton.

\ms \ni \underline{Case 2}. Choose $\C^6$-valued rational maps
$a_0, b_0, a_1, a_2$ on $\C$ such that \ben \item $\Im\pi_1$ is
spanned by $a_0, a_0', a_0'', a_0''', b_0$, \item $\Im\pi_2$ is
spanned by $a_0+ 2i\pi_1^\perp(a_1)$, $a_0'+ 2i\pi_1^\perp(a_1')$,
$a_0''+ 2i\pi_1^\perp(a_1'')$, \item $\Im\pi_3$ is spanned by $f_0
v_3 +f_1 D^{(1)}v_3$, where $v_3=\sum_{j=0}^2 (2i)^jP_{2j}(a_j)$,
$D^{(1)}v_3=\sum_{j=0}^2 (2i)^jP_{2j}(a_j')$, and $f_0,f_1: \C\to
\C$ are rational functions. \een
 These two cases give all $3$-unitons in $U(6)$ with rank data
 $(5,3,1)$. \eeg

\bs


\begin{thebibliography}{99}

\bibitem{A1997}
Anand, C.K., \emph{{W}ard's solitons},  Geom. Topol., {\textbf 1}
(1997), 9--20.

\bibitem{Ana98}
Anand, C.K., \emph{{W}ard's solitons II, Exact solutions}, Canad.
J. Math., {\textbf 50} (1998), 1119--1137.

\bibitem{BerGue91}
Bergvelt, M.J. and Guest, M., \emph{{A}ctions of loop groups on
harmonic maps}, Transactions AMS, {\textbf 326} (1991), 861--886.

\bibitem{BurGue97}
Burstall, F.E. and Guest, M.A., \emph{{H}armonic two-spheres in
compact symmetric spaces}, Math. Ann., {\textbf 309} (1997),
541--572

\bibitem{FokIoa98}
Fokas, A.S. and Ioannidou, T.A., \emph{{T}he inverse spectral
theory for the Ward equation and for the $2+1$ chiral model},
arXiv:hep-th/9806035.

\bibitem{I1996}
Ioannidou, T., \emph{{S}oliton solutions and nontrivial scattering
in an integrable chiral model in $(2+1)$ dimensions}, J. Math.
Phys., {\textbf 37} (1996), 3422--3441.
%
\bibitem{IW1995}
Ioannidou, T. and Ward, R.S., \emph{{C}onserved quantities for
integrable chiral model in $2+1$ dimensions}, Phys. Letters A,
{\textbf 208} (1995), 209--213.
%
\bibitem{IZ1998}
Ioannidou, T. and Zakrzewski, W., \emph{{S}olutions of the
modified chiral model in $(2+1)$ dimensions}, J. Math. Phys.,
{\textbf 39} (1998) no.5, 2693--2701.
%
\bibitem{TU2000}
Terng, C.L. and Uhlenbeck, K.,  \emph{{B}\"{a}cklund
transformations and loop group actions}, Comm. Pure Appl. Math.,
{\textbf 53} (2000), 1--75.
%
\bibitem{U1989}
Uhlenbeck, K., \emph{{H}armonic maps into Lie groups (classical
solutions of the chiral model)}, J. Diff. Geom., {\textbf 30}
(1989), 1--50.
%
\bibitem{V1990}
Villarroel, J., \emph{{T}he inverse problem for Ward's system},
Stud. Appl. Math., {\textbf 83} (1990), 211--222.
%
\bibitem{W1988}
Ward, R.S., \emph{{S}oliton solutions in an integrable chiral
model in $2+1$ dimensions}, J. Math. Phys., {\textbf 29} (1988),
386--389.
%
\bibitem{W1990}
Ward, R.S., \emph{{C}lassical solutions of the chiral model,
unitons, and holomorphic vector bundles}, Commun. Math. Phys.,
{\textbf 128} (1990), 319--332.
%
\bibitem{W1995}
Ward, R.S., \emph{{N}ontrivial scattering of localized solutions
in a $(2+1)$-dimensional integrable systems}, Phys. Letter A,
{\textbf 208} (1995), 203--208.
%
\bibitem{Woo1989}
Wood, J.C., \emph{{E}xplicit construction and parametrization of
harmonic maps into the unitary group}, Proc. London Math. Soc.,
{\textbf 58} (1989), 608--624.
%
\bibitem{Zh93}
Zhou, Z.X., \emph{{C}onstruction of explicit solutions of modified
principal chiral field in $1+2$ dimensions via Darboux
transformations},  Differential Geometry, edited by C.~H.~Gu et
al, pp.325--332, World Scientific, 1993.
%

\end{thebibliography}
\end{document}